\documentclass[12pt, reqno]{amsart}
\usepackage{amsmath, amsthm, amscd, amsfonts, amssymb, graphicx, color}
\usepackage[bookmarksnumbered, colorlinks, plainpages]{hyperref}
\input{mathrsfs.sty}
\usepackage[T1]{fontenc}
\usepackage[polish, english]{babel}
\usepackage[utf8]{inputenc}

\hypersetup{colorlinks=true,linkcolor=red, anchorcolor=green,
citecolor=cyan, urlcolor=red, filecolor=magenta, pdftoolbar=true}

\newtheorem{twr}{Theorem}[section]
\newtheorem{lem}[twr]{Lemma}

\newtheorem{cor}[twr]{Corollary}
\theoremstyle{definition}

\theoremstyle{remark}
\newtheorem{remark}[twr]{Remark}
\numberwithin{equation}{section}

\DeclareMathOperator{\absconv}{absconv}

\DeclareMathOperator{\lin}{lin}

\DeclareMathOperator{\tr}{Tr}
\DeclareMathOperator{\sgn}{sgn}
\DeclareMathOperator{\id}{Id}

\DeclareMathOperator{\rk}{rk}
\DeclareMathOperator{\re}{Re}
\DeclareMathOperator{\im}{Im}
\newcommand\dbm{d_{\mathrm{BM}}}
\usepackage{color}
\usepackage{enumerate}

\begin{document}

\setcounter{page}{1}

\title{Spaces with the maximal projection constant revisited}
\author[T. Kobos]
{Tomasz Kobos$^1$}
\address{$^1$Faculty of Mathematics and Computer Science, Jagiellonian University, \L ojasiewicza 6, 30-348 Krak\'ow, Poland}
\email{tomasz.kobos@uj.edu.pl}

\begin{abstract}
Let $n \geq 2$ be an integer such that an equiangular set of vectors $w_1, \ldots, w_d$ of the maximal possible cardinality (that is, attaining the classical Gerzon upper bound) exists in $\mathbb{K}^n$, where $\mathbb{K}=\mathbb{R}$ or $\mathbb{K}=\mathbb{C}$ (so that $d=\frac{n(n+1)}{2}$ in the real case and $d=n^2$ in the complex case). We provide a complete characterization of $n$-dimensional normed spaces whose absolute projection constant is maximal among all $n$-dimensional normed spaces over $\mathbb{K}$. The characterization states that $X$ has the maximal projection constant if and only if it is isometric to a space whose dual unit ball is contained between the absolutely convex hull of the vectors $w_1, \ldots, w_d$ and a suitably rescaled zonotope generated by the same vectors. As a consequence, we obtain that, in the considered situations, $n=2$ with $\mathbb{K}=\mathbb{R}$ is the only case in which there is, up to isometry, a unique norm on $\mathbb{K}^n$ with the maximal projection constant. In this case, the unit ball is a linear image of a regular hexagon in $\mathbb{R}^2$.
\end{abstract} 
\maketitle

\section{Introduction}

Let $X$ be a normed space over $\mathbb{K}$ (where $\mathbb{K}=\mathbb{R}$ or $\mathbb{K}=\mathbb{C}$) and let $Y \subseteq X$ be its linear subspace. By $\lambda(Y, X)$ we denote the \emph{relative projection constant} of $Y$, defined as the infimum of $\|P\|$, where $P: X \to Y$ is a linear projection. If a linear projection $P: X \to Y$ satisfies $\|P\|=\lambda(Y, X)$, it is called a \emph{minimal projection}. By $\lambda(X)$ we shall denote the \emph{absolute projection constant} of $X$, which is defined as the supremum of $\lambda(T(X), Z)$, where $Z$ is any normed space such that there exists a linear isometric embedding $T: X \to Z$. We shall consider only the finite-dimensional setting. In this case, it follows from a standard compactness argument that for any $Y \subseteq X$ there exists a minimal projection onto $Y$. Moreover, since by a famous result of Kadec and Snobar \cite{kadecsnobar} we always have $\lambda(X) \leq \sqrt{n}$ for $\dim X = n$, it is also easy to deduce that the maximal possible absolute projection constant $\lambda(X)$ among all $n$-dimensional normed spaces over $\mathbb{K}$ is attained by some $n$-dimensional normed space. We shall denote this maximal projection constant by $\lambda_{\mathbb{K}}(n)$. Since $\lambda(X)=1$ holds for any $1$-dimensional space $X$, we have $\lambda_{\mathbb{K}}(1)=1$ and throughout the paper we will therefore always assume $n \geq 2$. It is well-known that for every finite-dimensional normed space $X$ there exists an isometric embedding $T: X \to \ell_{\infty}$ and for every such isometric embedding we have $\lambda(T(X), \ell_{\infty})=\lambda(X)$. This means that the maximal projection constant spaces can be found among subspaces of $\ell_{\infty}$ of a given dimension. If the norm of $X$ is given by $\max_{ 1 \leq i \leq N} |f_i(x)|$ for $x \in \mathbb{K}^n$, where $f_i: \mathbb{K}^n \to \mathbb{K}$ are linear functionals and $N \geq n$ is an integer, then $X$ can be isometrically embedded into $\ell_{\infty}^N$ and again, the absolute projection constant of $X$ is the relative projection constant of the image when considered as a subspace of $\ell_{\infty}^N$. In the real case, this happens precisely when the unit ball of $X$ is a convex polytope, and therefore we shall refer to such norms on $\mathbb{K}^n$ as \emph{polyhedral norms}, both in the real and complex settings. It should be noted, however, that in sharp contrast to the real case, complex polyhedral norms not only fail to form a self-dual class, but the dual norm of a polyhedral norm on $\mathbb{C}^n$ (for $n \geq 2$) is never polyhedral (see Section~2 in \cite{koboslewicki}).

Projection constants of normed spaces have been extensively studied for several decades due to their deep connections with other fundamental concepts in the isometric theory of Banach spaces. For a general overview, we refer the reader to Section~$8$ in \cite{tomczak} or Section~III.B in \cite{wojtaszczyk}. While projection constants have been investigated in various contexts, this paper focuses on determining the maximal absolute projection constant among all spaces of a given dimension. This line of research has a rich history of remarkable results, as well as several incorrect proofs published over the years. Our main goal is to provide a complete characterization of $n$-dimensional spaces with the maximal projection constant in all situations where the exact value of $\lambda_{\mathbb{K}}(n)$ is currently known. To contextualize our investigation, we provide a historical overview of previous research. To clarify the current state of knowledge, it is also necessary to address certain arguments from the literature that ultimately proved to be incorrect. We restrict our summary below to the papers most directly relevant to our work.

We start our overview with a 1960 paper of Gr\"unbaum ``Projection Constants'' (\cite{grunbaum}), where some of the most fundamental facts about projection constants were established and where the question of determining the maximal projection constant of spaces of a given dimension was first raised. Among other things, Gr\"unbaum calculated absolute projection constants of some two-dimensional real normed spaces and observed that the projection constant of a two-dimensional normed space $X$ with a regular hexagon as the unit ball satisfies $\lambda(X) = \frac{4}{3}$. He conjectured that this is the exact value of $\lambda_{\mathbb{R}}(2)$. This turned out to be a highly challenging problem, as it was not until more than three decades later that it was first tackled by K\"onig and Tomczak-Jaegermann in their 1994 paper ``Norms of Minimal Projection'' \cite{konigtomczak}. The authors claimed a much more general result (Theorem $1.1$ (a)), that for every dimension $n \geq 1$ we have $\lambda_{\mathbb{K}}(n) \leq \delta_{\mathbb{K}}(n)$, where
$$\delta_{\mathbb{K}}(n)= \begin{cases} \frac{2}{n+1} \left ( 1 + \frac{n-1}{2} \sqrt{n+2}\right ) \quad \hbox{ for } \mathbb{K}=\mathbb{R}, \\
\frac{1}{n} \left ( 1 + (n-1) \sqrt{n+1} \right ) \quad \hbox{ for } \mathbb{K}=\mathbb{C}.
\end{cases}.$$
This estimate includes the Gr\"unbaum's conjecture as a particular case for $n=2$ and $\mathbb{K}=\mathbb{R}$. Moreover, this paper was likely the first to establish a connection between projection constants and equiangular sets of vectors. We recall that a set of pairwise distinct unit vectors $w_1, \ldots, w_m \in \mathbb{K}^n$ is called \emph{equiangular} if the value $|\langle w_i, w_j \rangle|$ is constant for all $1 \leq i < j \leq m$. The study of equiangular sets of vectors and equiangular tight frames is of independent interest and has gained significant attention in recent decades. It is a classical fact attributed to Gerzon (see \cite{lemmensseidel}) that the cardinality of an equiangular set of vectors in $\mathbb{K}^n$ is less than or equal to $d_{\mathbb{K}}(n)$, where
$$d_{\mathbb{K}}(n)= \begin{cases} \frac{n(n+1)}{2} \quad \hbox{ for } \mathbb{K}=\mathbb{R}, \\
n^2 \quad \hbox{ for } \mathbb{K}=\mathbb{C}.
\end{cases}.$$
In the real case, this upper bound is achieved for $n=2, 3, 7, 23$, but no other examples are currently known. Available evidence suggests that these may be the only dimensions where equality holds in the real setting. The complex case, however, is vastly different. The upper bound $n^2$ is currently known to be attained in more than 20 different complex dimensions, and a famous conjecture due to Zauner \cite{zauner} asserts that it is optimal in every dimension. Zauner's conjecture remains a major open problem in the field of quantum information theory, where it is frequently stated in terms of the existence of a SIC-POVM (a symmetric, informationally complete, positive operator-valued measure). In part (b) of Theorem 1.1 from \cite{konigtomczak}, K\"onig and Tomczak-Jaegermann stated that the upper bound $\lambda_{\mathbb{K}}(n) \leq \delta_{\mathbb{K}}(n)$ is attained if and only if there exists an equiangular set of vectors in $\mathbb{K}^n$ of maximal cardinality $d_{\mathbb{K}}(n)$. Furthermore, they claimed that any space with the maximal projection constant has the dual unit ball equal to the absolutely convex hull of such a maximal equiangular set. Recall that the \emph{absolutely convex hull}of a set $C \subseteq \mathbb{K}^n$ is defined similarly for real and complex scalars as
$$\absconv C = \left\{ \sum_{j=1}^{N} \alpha_j x_j : N \geq 1, \, x_j \in C, \, \alpha_j \in \mathbb{K}, \, \sum_{j=1}^{N} \vert{}\alpha_j\vert{} \le 1 \right\}.$$
Part (c) further stated that if $\mathbb{K}=\mathbb{R}$ and $n \in \{2, 3, 7, 23\}$, then there is a unique (up to isometry) $n$-dimensional normed space with the maximal projection constant -- namely, the space whose dual unit ball is the absolutely convex hull of a maximal equiangular set of vectors in $\mathbb{R}^n$ (since in these cases, a maximal equiangular set of vectors is unique up to an orthogonal transformation of $\mathbb{R}^n$).

Thus, it would seem that Theorem~1.1 from \cite{konigtomczak} completely settled the problem of determining and realizing maximal projection constants in all spaces $\mathbb{K}^n$ for which an equiangular set of cardinality $d_{\mathbb{K}}(n)$ exists. However, it later turned out that the argument was incorrect and, interestingly, the mistake was twofold. The crucial step in the proof of Theorem~1.1 was Proposition~3.1, which, without going into too much detail, consisted of two separate parts: one dealing with the non-zero case ($\mu_s^0 \neq 0$) and one dealing with the zero case ($\mu_s^0 = 0$). The zero case was handled in just a few lines immediately following the statement of Proposition~3.1, whereas the proof of the non-zero case relied on a highly technical argument involving Lagrange multipliers.

Both parts of this proposition eventually turned out to be incorrect. The failure of the short argument for the zero case was discovered by the authors themselves, who, nine years later, published the paper ``Spaces with maximal projection constants'' \cite{konigtomczak2}. The results of this paper directly contradict part~(c) of Theorem~1.1, which stated that for $n \in \{2, 3, 7, 23\}$ there is a unique $n$-dimensional real normed space with the maximal projection constant. In contrast, the authors stated in \cite{konigtomczak2} that although in $\mathbb{R}^2$ the normed space with a regular hexagon as its unit ball is the unique normed space with the maximal projection constant, there exist two non-isometric norms on $\mathbb{R}^3$ and on $\mathbb{C}^2$ attaining the maximal projection constant. It is worth noting that it is not entirely clear on what grounds uniqueness was claimed for $\mathbb{R}^2$, as no proof of this assertion appears in the paper. The authors explicitly acknowledged the flaw concerning the zero case in the proof of Proposition~5 (page 362). This mistake was responsible for the non-uniqueness phenomenon, whereas part~(a) of Theorem~1.1 from \cite{konigtomczak} appeared to remain intact. It is important to emphasize that the results of the second paper \cite{konigtomczak2} still relied heavily on those of the earlier work. In particular, Proposition~5, which is crucial for the arguments in \cite{konigtomczak2}, was proved using the non-zero case of Proposition~3.1 from \cite{konigtomczak}.

A further important development occurred in 2009, when Chalmers and Lewicki published the paper ``Three-dimensional subspace of $\ell_{\infty}^5$ with maximal projection constant'' \cite{chalmerslewicki}. In this work, the authors showed that the non-zero part of Proposition~3.1 from \cite{konigtomczak} is also incorrect. The computational error was identified at the end of page 554. Moreover, the main result of the paper showed that not only was the proof of Proposition~3.1 incorrect, but its statement was in fact false in general (Remark~3.1 on page 591). Consequently, this also invalidated the results from \cite{konigtomczak2}, which relied indirectly on Proposition~3.1. In a subsequent paper \cite{chalmerslewicki2}, Chalmers and Lewicki provided a new proof of the Gr\"unbaum conjecture, thereby settling the problem. However, their proof was long and highly technical. It did not appear amenable to adaptation either to establish part~(a) of Theorem~1.1 from \cite{konigtomczak}, namely the estimate $\lambda_{\mathbb{K}}(n)\leq \delta_{\mathbb{K}}(n)$ for all $n$, or to determine all norms in $\mathbb{R}^2$ with the maximal projection constant. Thus, at that time, $\mathbb{R}^2$ remained the only non-trivial case for which the maximal absolute projection constant had been determined. Moreover, although a norm whose unit ball is a regular hexagon was known to attain the maximal projection constant in $\mathbb{R}^2$, it was not known whether this norm was unique (up to isometry).

In 2019, Basso published the paper ``Computation of maximal projection constants'' \cite{basso}, in which he derived a new formula for the maximal projection constant. This enabled him to obtain an alternative proof of the Gr\"unbaum conjecture based on a graph-theoretic argument. However, it still did not seem feasible to extract the equality case from this approach. Furthermore, Basso claimed in Theorem~1.4 of \cite{basso} that for every $n \geq 1$ there exists an $n$-dimensional real polyhedral normed space attaining the maximal projection constant. Equivalently, the theorem asserted that the maximal projection constant in $\mathbb{R}^n$ can always be attained by an $n$-dimensional subspace of $\ell_{\infty}^N$ for sufficiently large $N$. If true, such a result would have important consequences for the study of maximal projection constants, since polyhedral norms are often considerably easier to analyze than arbitrary norms. Unfortunately, this statement also turned out to be false (see Remark~\ref{remprojconst} for further details).

A major breakthrough came in 2023, when Deręgowska and Lewandowska published the paper ``A simple proof of the Gr\"unbaum conjecture'' \cite{deregowskalewandowska}, which led to significant progress on long-standing problems concerning maximal projection constants. In this work, the authors observed that a result of Bukh and Cox \cite{bukhcox} on isotropic measures can be adapted to provide a proof of the inequality $\lambda_{\mathbb{K}}(n)\leq \delta_{\mathbb{K}}(n)$ for all $n$, thereby establishing part~(a) of Theorem~1.1 from \cite{konigtomczak}. The resulting proof is remarkably short and simple, with the main estimate relying only on the Cauchy--Schwarz inequality and basic linear algebra. This leaves the remaining parts of Theorem~1.1 from \cite{konigtomczak}, which concern the characterization of the equality case and the uniqueness of spaces with the maximal projection constant. As we shall see, our main result will in particular establish that equality in the estimate $\lambda_{\mathbb{K}}(n)\leq \delta_{\mathbb{K}}(n)$ holds precisely when there exists an equiangular set of vectors in $\mathbb{K}^n$ of cardinality $d_{\mathbb{K}}(n)$. The problems of characterization and uniqueness will then be resolved in Theorem~\ref{twrchar} and Corollary~\ref{coruniq}. We note that even though part~(a) of Theorem~1.1 from \cite{konigtomczak} has now been proved correctly, this does not automatically validate the results of \cite{konigtomczak2}, since they depend heavily on the incorrect Proposition~3.1. The main goal of the present paper is not only to settle the question of uniqueness, but also to provide a complete classification of all spaces with the maximal projection constant in every dimension for which maximal equiangular sets are known to exist -- precisely the dimensions where the value of $\lambda_{\mathbb{K}}(n)$ is currently known. For instance, even for $\mathbb{R}^4$, determining the exact maximal projection constant remains an open problem.

It is well known (see Theorems~2.1 and~2.2 in \cite{chalmerslewicki} or Appendix~A in \cite{foucartskrzypek}) that the maximal absolute projection constant can be expressed as
\begin{equation} \label{formula} \lambda_{\mathbb K}(n) = \sup \sum_{i,j=1}^{N} t_it_j|\langle u_i,u_j\rangle|, 
\end{equation} 
where the supremum is taken over all integers $N\ge1$, all non-negative real numbers $t_i$ satisfying $\sum_{i=1}^{N} t_i^2=1$, and all vectors $u_1,\ldots,u_N\in\mathbb K^n$ forming a tight frame with constant $1$ in $\mathbb K^n$. Recall that vectors $u_1,\ldots,u_N\in\mathbb K^n$ form a \emph{tight frame with constant $\alpha$} if, for every $x\in\mathbb K^n$, 
$$ x=\alpha\sum_{i=1}^{N}\langle x,u_i\rangle u_i. $$
A particularly important special case is when the vectors $u_i$ are unit vectors forming an equiangular set. In this situation, such a tight frame is called an \emph{equiangular tight frame}, or simply an \emph{ETF}. By the previously mentioned upper bound on the cardinality of an equiangular set of vectors, it follows that if $u_1,\ldots,u_N\in\mathbb K^n$ form an ETF, then $N\le d_{\mathbb K}(n)$.
In the case $N=d_{\mathbb{K}}(n)$, we shall say that an ETF is a \emph{maximal ETF}. It is known (see Theorem~5.10 in \cite{foucart}) that every equiangular set of vectors $u_1,\ldots,u_N$ in $\mathbb{K}^n$ of cardinality $N=d_{\mathbb{K}}(n)$ is automatically an ETF. Moreover, for all $1 \leq i < j \leq N$ we have $|\langle u_i, u_j \rangle| = \varphi_{\mathbb{K}}(n)$ (see, for example, Theorem $5.7$ in \cite{foucart}), where
$$\varphi_{\mathbb{K}}(n) = \begin{cases} \frac{1}{\sqrt{n+2}} \quad \hbox{ for } \mathbb{K}=\mathbb{R}, \\
\frac{1}{\sqrt{n+1}} \quad \hbox{ for } \mathbb{K}=\mathbb{C}.
\end{cases}$$

To determine the value of $\lambda_{\mathbb{K}}(n)$, it suffices to estimate the expression appearing on the right-hand side of \eqref{formula} over all integers $N \geq 1$. However, a characterization of all spaces attaining $\lambda_{\mathbb{K}}(n)$ requires more than an analysis of the equality case for each fixed $N$, as this would only yield a description of the maximizers within the class of polyhedral norms. To obtain an upper bound for the right-hand side of \eqref{formula}, Deręgowska and Lewandowska used an approximation argument that allowed them to assume that all weights are equal, i.e. $t_i = \frac{1}{\sqrt{N}}$ for every $1 \leq i \leq N$ (see Remark \ref{remprojconst} for details). Moreover, the zero vectors $u_i$ were easily discarded in the proof, as they do not contribute positively to the right-hand side. This allowed for a streamlined version of a more general estimate due to Bukh and Cox \cite{bukhcox} for isotropic measures (see Lemma~5 in \cite{bukhcox}). However, neither of these simplifications can be used when characterizing the equality case. In the result below, we provide an estimate for the right-hand side of \eqref{formula} involving the weights $t_i$ and, most importantly, we give a complete characterization of the equality case. Here and throughout the paper, by $\|x\|=\sqrt{\sum_{i=1}^{n} |x_i|^2}$ we shall always denote the Euclidean norm of a vector $x \in \mathbb{K}^n$ and by a \emph{unimodular scalar} we mean a scalar $z \in \mathbb{K}$ with $|z|=1$.

\begin{twr}
\label{twrbound}
Let $n \geq 2$ be an integer. Suppose that $N \geq n$ is an integer and let $u_1, \ldots, u_N \in \mathbb{K}^n$ be vectors forming a tight frame in $\mathbb{K}^n$ with constant $1$. Let $t \in \mathbb{R}^N$ be a vector with non-negative coordinates satisfying $\|t\|=1$. Then the inequality
$$\sum_{i, j=1}^{N} t_it_j |\langle u_i, u_j \rangle| \leq \delta_{\mathbb{K}}(n)$$
is true. Moreover, equality holds if and only if there exists a maximal ETF $w_1, \ldots, w_{d}$ in $\mathbb{K}^n$ (where $d=d_{\mathbb{K}}(n)$) and a partition $\{1, \ldots, N\} = A_0 \sqcup A_1 \sqcup \ldots \sqcup A_d$ satisfying the following conditions:
\begin{enumerate}
\item $A_0 = \{1 \leq i \leq N: u_i=0\}$ and for all $1 \leq i \leq N$ we have $t_i=0$ if and only if $i \in A_0$.
\item For every $1 \leq j \leq d$ and every $i \in A_j$ there exists a unimodular scalar $\alpha_i$ such that 
$$\frac{u_i}{\|u_i\|} = \alpha_i w_j.$$
\item For $1 \leq j \leq d$ we have
$$\sum_{i \in A_j} t_i^2 = \frac{1}{d}$$
(in particular we must have $N \geq d$, as each set $A_j$ has to be non-empty).
\item $\|u_i\|=\sqrt{n}t_i$ for every $1 \leq i \leq N$.
\end{enumerate}
\end{twr}
The result implies that in the equality case, a maximal ETF must exist in $\mathbb{K}^n$. We emphasize that the estimate itself is not new. The proof is carried out in a more general setting. This additional generality is not essential if one is interested only in the estimate itself, but it becomes crucial for the characterization of spaces with the maximal projection constant, where the equality conditions for the general case play a fundamental role. As expected, our reasoning follows very closely that of \cite{deregowskalewandowska}. Some adjustments of certain steps are necessary to accommodate a general vector $t$ and the possibility that some vectors $u_i$ are equal to $0$, but the essential idea remains the same. A detailed analysis of the estimates appearing in the proof yields the necessary conditions for equality, which are then easily shown to be sufficient as well. The proof of Theorem~\ref{twrbound} is presented in Section~\ref{seceqcond}.

Part (4) of the equality conditions given above corresponds to the non-zero case of the generally incorrect Proposition~3.1 from \cite{konigtomczak} and shows that its conclusion remains valid in certain situations. One could argue that this observation alone might be sufficient to recover some of the results from \cite{konigtomczak2}, at least in some special cases. However, a great deal of caution is needed here, given that Proposition~3.1 is incorrect in full generality. In order to answer the question of uniqueness, we build the argument from scratch, especially since our main goal is to provide a complete characterization of spaces with the maximal projection constant (in the situations where a maximal ETF exists), which was not attempted in \cite{konigtomczak2}. In the characterization below, for vectors $x_1,\ldots,x_N\in\mathbb K^n$, we denote by $Z(x_1,\ldots,x_N)$ the Minkowski sum of the sets $\absconv\{x_i\}$, that is,
$$Z(x_1, \ldots, x_N) = \left \{a_1x_1 + \ldots + a_Nx_N: \ a_i \in \mathbb{K}, \ |a_i| \in [0, 1] \ \text{ for } \ 1 \leq i \leq N \right \}.$$
In the real case, the set $Z(x_1,\ldots,x_N)$ is a centrally symmetric zonotope (and hence, in particular, a convex polytope). Our characterization is as follows.
\begin{twr}
\label{twrchar}
Let $n \geq 2$ be an integer and let $X = (\mathbb{K}^n, \| \cdot \|)$ be a normed space. Then $X$ satisfies $\lambda(X) = \delta_{\mathbb{K}}(n)$ if and only if there exists a maximal ETF $w_1, \ldots, w_d \in \mathbb{K}^n$ and a linear transformation $T: \mathbb{K}^n \to \mathbb{K}^n$ such that
$$\absconv \{w_1, \ldots, w_d\} \subseteq T(B_{X^*}) \subseteq \frac{n}{d\delta_{\mathbb{K}}(n)} Z(w_1, \ldots, w_d).$$
\end{twr}

Let us note that the theorem implies that if there exists an $n$-dimensional space with $\lambda(X) = \delta_{\mathbb{K}}(n)$, then a maximal ETF must exist in $\mathbb{K}^n$, and furthermore the inclusion 
\begin{equation}
\label{inclusion}
\absconv \{w_1, \ldots, w_d\} \subseteq \frac{n}{d\delta_{\mathbb{K}}(n)} Z(w_1, \ldots, w_d)
\end{equation}
holds. In Lemma~\ref{lemincl} we shall prove that the case $\mathbb R^2$ is the only situation in which this inclusion is not strict. This is exactly the reason why uniqueness of the space with the maximal projection constant holds in $\mathbb R^2$, but not in the remaining cases (see Corollary~\ref{coruniq}).

The norms on $\mathbb{K}^n$ dual to those with the unit balls $\absconv \{w_1, \ldots, w_d\}$ and $\frac{n}{d\delta_{\mathbb{K}}(n)} Z(w_1, \ldots, w_d)$ appeared explicitly in \cite{konigtomczak2}. They were denoted by $X_n$ and $Y_n$, respectively, although their definition was more general, as the authors attempted to cover all dimensions $n$. It was claimed that both spaces attain the maximal projection constant while being, in general, non-isometric. However, to the best of our understanding, a full characterization such as that given in Theorem~\ref{twrchar} was neither formulated nor attempted. These spaces were described explicitly in the case of $\mathbb{R}^3$ in Proposition~2, and it was stated that all spaces between them have maximal projection constant. The authors do not appear to have claimed, however, that these are all possible spaces in $\mathbb{R}^3$ with the maximal projection constant. Furthermore, the relation between these two norms (corresponding to the inclusion~\eqref{inclusion}) does not appear to have been observed in general in situations where a maximal ETF exists. The non-uniqueness of spaces with the maximal projection constant was stated in Theorem~3, although the range of dimensions considered there misses numerous small cases (see the discussion before Corollary~\ref{coruniq} for details).

The proof of Theorem~\ref{twrchar} combines ideas from \cite{konigtomczak} and \cite{konigtomczak2}, building the argument from scratch and, in the process, reestablishing the estimate \eqref{formula}. It should be emphasized that, although the results of \cite{konigtomczak} and \cite{konigtomczak2} cannot be regarded as correct, many of the underlying ideas play an important role in our approach. The proof of Theorem~\ref{twrchar} is spread across Sections~\ref{seccharpolyhedral} and \ref{secchargeneral}. In Section~\ref{seccharpolyhedral} we first establish the characterization in the case where the norm of $X$ is polyhedral. This situation is easier to handle, since in this case $X$ can be isometrically embedded into $\ell_{\infty}^N$ for a suitable integer $N \geq 1$. In Section~\ref{secchargeneral} we extend the characterization to the general case by approximating an arbitrary norm on $\mathbb{K}^n$ with a sequence of polyhedral norms. While it is straightforward to deduce from the polyhedral case that a space $X$ satisfying the inclusions from Theorem~\ref{twrchar} has the maximal projection constant, the proof of the converse implication requires considerably more effort. This is because it necessitates a refinement of most of the arguments used in the polyhedral case in order to adapt them to the situation where a general space $X$ with the maximal projection constant is approximated by polyhedral spaces with almost maximal projection constants. This additional work is necessary to show that the characterization established in the polyhedral setting carries over to the limit.

Throughout the paper, for a non-zero $z \in \mathbb{K}$, we denote by $\sgn(z)$ the unique unimodular scalar satisfying $z\sgn(z) = |z|$. Additionally, we define $\sgn(0)=1$. We shall frequently use the following simple observation: if vectors $u_1,\ldots,u_N \in \mathbb{K}^n$ form a tight frame with constant $1$, then
\begin{equation}
\label{eqtrace}
\sum_{i=1}^{N} \|u_i\|^2 = n.
\end{equation}
This follows immediately by comparing the traces of the operators appearing in the equality $x = \sum_{i=1}^{N} \langle x, u_i \rangle u_i$. Moreover, taking the inner product with $x$ in the same equality yields
\begin{equation}
\label{eqnorma}
\|x\|^2 = \sum_{i=1}^{N} |\langle x, u_i \rangle |^2,
\end{equation}
for every $x \in \mathbb{K}^n$.

For a matrix $A$ with entries in $\mathbb{K}$, we denote by $A^*$ the transpose when $\mathbb{K}=\mathbb{R}$ and the Hermitian transpose when $\mathbb{K}=\mathbb{C}$. The same notation is used for the adjoint of a linear operator over $\mathbb{K}$. By $e_i$ we denote the $i$-th canonical basis vector, and by $\langle \cdot,\cdot\rangle$ the usual inner product, assumed to be linear in the first variable and conjugate-linear in the second. Finally, for an integer $N \geq 1$ by $\ell_2^N$ and $\ell_{\infty}^N$ denote the $N$-dimensional normed spaces over $\mathbb{K}$ equipped with the Euclidean norm and the maximum norm, respectively.

\section{Equality conditions for the expression bounding the absolute projection constant}
\label{seceqcond}

\begin{proof}[Proof of Theorem \ref{twrbound}]
We begin by proving the estimate and establishing the necessity of the equality conditions. For the sake of clarity, let us first assume that $A_0=\emptyset$. We do not, however, assume that all coordinates of $t$ are non-zero. For each $1 \leq i \leq N$ we define a linear operator $L_i: \mathbb{K}^n \to \mathbb{K}^n$ of rank one as
$$L_i(x) = \sqrt{t_i} \|u_i\|^{\frac{-3}{2}} \langle x, u_i \rangle u_i,$$
i.e. $L_i = \sqrt{t_i} \|u_i\|^{\frac{-3}{2}} u_iu_i^{*}$ in the matrix notation. Let $\mathcal{H} = \lin \{\id_n, L_1, \ldots L_N \}$. We note that when $\mathbb K=\mathbb R$, the matrix $L_i$ is symmetric and hence $\mathcal{H}$ is a linear subspace of $d_{\mathbb{R}}(n)=\frac{n(n+1)}{2}$-dimensional vector space of real symmetric matrices of size $n \times n$. In the case of $\mathbb{K}=\mathbb{C}$ it is a Hermitian matrix, but since the Hermitian matrices do not form a linear subspace over $\mathbb{C}$, we consider $\mathcal{H}$ as a subspace of $d_{\mathbb{C}}(n)=n^2$-dimensional vector space of complex matrices of size $n \times n$. In either case, $\mathcal H$ can be regarded as a Hilbert space over $\mathbb{K}$ with the Frobenius inner product given as $\langle A, B \rangle_F = \tr(AB^*)$ for $A, B \in \mathcal{H}$. We observe that for any $1 \leq i, j \leq N$ we have
$$\langle L_i, L_j \rangle_F = \tr(L_iL^*_j) = \tr (L_iL_j) = \tr (\sqrt{t_it_j} \|u_i\|^{\frac{-3}{2}} \|u_j\|^{\frac{-3}{2}} u_iu_i^*u_ju_j^*)$$
$$=\sqrt{t_it_j}\|u_i\|^{\frac{-3}{2}} \|u_j\|^{\frac{-3}{2}} |\langle u_i, u_j \rangle |^2.$$
We now define an operator $G:\ell_2^N\to\ell_2^N$ by
$$G(x) = \sum_{j=1}^{N} \left ( \sum_{i=1}^{N} x_i \langle L_i, L_j \rangle_F \right ) e_j.$$
Thus, for any $1 \leq i, j \leq N$ we have $\langle G(e_i), e_j \rangle = \langle L_i, L_j \rangle_F$ and the operator $G$ may be regarded as the Gram matrix of vectors $L_1, \ldots, L_N$ with respect to the Frobenius inner product. We shall prove that $G$ is of rank at most $\dim \mathcal{H} \leq d_{\mathbb{K}}(n)$. To this end, let us observe that there is a decomposition of the form $G = SS^*$, where a linear operator $S: \mathcal{H} \to \ell_2^N$ is defined as
$$S(L) = \sum_{i=1}^{N} \langle L, L_i \rangle_F e_i,$$
for $L \in \mathcal{H}$.
Thus, for the adjoint operator $S^*: \ell_2^N \to \mathcal{H}$ we have $S^*(x) = \sum_{i=1}^{N} x_i L_i$ for $x \in \ell_2^N$. It is now simple to check algebraically that indeed $G=SS^*$. In particular, the image of $G$ is contained in a subspace $V=S(\mathcal{H})$ whose dimension does not exceed that of $\mathcal H$. We further observe that if we define a linear operator $T: \ell_2^N \to \ell_2^N$ as
$$T(x) = \sum_{j=1}^{N} \left ( \sum_{i=1}^{N} x_i \sqrt{t_it_j \|u_i\| \|u_j\|} \right ) e_j,$$
or equivalently, $\langle T(e_i), e_j \rangle = \sqrt{t_it_j \|u_i\| \|u_j\|}$ for $1 \leq i, j \leq N$, then the image of $T$ is also contained in a subspace $V$. Indeed, it is enough to prove that $T(e_j) \in V$ for every $1 \leq j \leq N$. However, we note that for $L = \sqrt{t_j \|u_j\|} \id_n \in \mathcal{H}$ we have
$$S(L) = \sum_{i=1}^{N} \langle L, L_i \rangle_F e_i = \sqrt{t_j \|u_j\|} \sum_{i=1}^{N}  \left ( \sqrt{t_i} \|u_i\|^{\frac{-3}{2}} \tr(u_iu_i^*) \right ) e_i$$
$$=\sqrt{t_j \|u_j\|} \sum_{i=1}^{N} \left (  \sqrt{t_i} \|u_i\|^{\frac{-3}{2}} \|u_i\|^2 \right ) e_i = \sum_{i=1}^{N} \sqrt{t_i t_j \|u_i\| \|u_j\| } e_i = T(e_j).$$
This proves that $T(e_j) \in V$ for every $1 \leq j \leq N$ and consequently $T(\ell_2^N) \subseteq V$. In particular, for any $a, b \in \mathbb{R}$ a self-adjoint operator $A: \ell_2^N \to \ell_2^N$ defined as $A = a G - b T$ satisfies $A(\ell_2^N) \subseteq V$. Hence, such an operator is of rank at most $\dim V \leq \dim \mathcal{H} \leq d$ and therefore
\begin{equation}
\label{trdim}
\tr(A^2)=\tr(A^*A) \geq \frac{\tr(A)^2}{\rk A} \geq \frac{\tr(A)^2}{\dim \mathcal{H}} \geq \frac{\tr(A)^2}{d}.
\end{equation}
Expressing the traces in the canonical basis, we obtain
$$\sum_{i, j=1}^{N} \left ( a \sqrt{t_it_j}\|u_i\|^{\frac{-3}{2}} \|u_j\|^{\frac{-3}{2}} |\langle u_i, u_j \rangle |^2 - b \sqrt{t_it_j \|u_i\| \|u_j\|}  \right )^2 \geq \frac{\left (  \sum_{i=1}^{N} (a-b) t_i \|u_i\|   \right )^2}{d}$$
or
\begin{equation}
\label{traceestimate}
\sum_{i, j=1}^{N} t_i t_j \frac{  \left (a |\langle u_i, u_j \rangle|^2 - b \|u_i\|^2 \|u_j\|^2 \right )^2}{\|u_i\|^3\|u_j\|^3} \geq  \frac{(a-b)^2 \left (  \sum_{i=1}^{N}  t_i \|u_i\|   \right )^2}{d}.
\end{equation}
For simplicity, let us write $\varphi = \varphi_{\mathbb{K}}(n)$. We shall focus on the expression
$$\sum_{i, j=1}^{N} t_i t_j \frac{\left (|\langle u_i, u_j \rangle| - \varphi \|u_i\| \|u_j\| \right )^2}{\|u_i\| \|u_j\|}.$$
We observe that
$$\sum_{i, j=1}^{N} t_i t_j \frac{\left (|\langle u_i, u_j \rangle| - \varphi \|u_i\| \|u_j\| \right )^2}{\|u_i\| \|u_j\|}=\sum_{i, j=1}^{N} t_it_j \frac{\left (|\langle u_i, u_j \rangle|^2 - \varphi^2 \|u_i\|^2 \|u_j\|^2 \right )^2}{\|u_i\| \|u_j\|\left (|\langle u_i, u_j \rangle| + \varphi  \|u_i\| \|u_j\| \right )^2}.$$
By the Cauchy--Schwarz inequality we have $| \langle u_i, u_j \rangle | \leq \|u_i\| \|u_j\|$ for every $1 \leq i, j \leq N$. Hence, we can estimate
$$
\sum_{i, j=1}^{N} t_it_j \frac{\left (|\langle u_i, u_j \rangle|^2 - \varphi^2 \|u_i\|^2 \|u_j\|^2 \right )^2}{\|u_i\| \|u_j\|\left (|\langle u_i, u_j \rangle| + \varphi  \|u_i\| \|u_j\| \right )^2} \geq \sum_{i, j=1}^{N} t_it_j \frac{\left (|\langle u_i, u_j \rangle|^2 - \varphi^2 \|u_i\|^2 \|u_j\|^2 \right )^2}{(1+\varphi)^2\|u_i\|^3\|u_j\|^3}.
$$
If equality holds in this estimate, then for each summand either equality holds in the Cauchy--Schwarz inequality or the numerator vanishes. Hence, if equality holds, then for all $1 \leq i, j \leq N$ we have
\begin{equation}
\label{condequi}
\left | \left \langle \frac{u_i}{\|u_i\|}, \frac{u_j}{\|u_j\|} \right \rangle \right | \in \{1, \varphi\}.
\end{equation}
Applying \eqref{traceestimate} with $a=\frac{1}{1+\varphi}$ and $b = \frac{\varphi^2}{1+\varphi}$ we can further estimate
\begin{align*}
\sum_{i, j=1}^{N} t_it_j \frac{\left (|\langle u_i, u_j \rangle|^2 - \varphi^2 \|u_i\|^2 \|u_j\|^2 \right )^2}{(1+\varphi)^2\|u_i\|^3\|u_j\|^3} &= \sum_{i, j=1}^{N} t_i t_j \frac{  \left (\frac{|\langle u_i, u_j \rangle|^2}{1+\varphi} - \frac{\varphi^2}{1+\varphi} \|u_i\|^2 \|u_j\|^2 \right )^2}{\|u_i\|^3\|u_j\|^3} \\ 
&\geq \frac{(1 - \varphi)^2 \left (  \sum_{i=1}^{N}  t_i \|u_i\|   \right )^2}{d}.
\end{align*}
To summarize, we have proved that
\begin{equation}
\label{mainestimate}
\sum_{i, j=1}^{N} t_i t_j \frac{\left (|\langle u_i, u_j \rangle| - \varphi \|u_i\| \|u_j\| \right )^2}{\|u_i\| \|u_j\|} \geq\frac{(1 - \varphi)^2 \left (  \sum_{i=1}^{N}  t_i \|u_i\|   \right )^2}{d}.
\end{equation}
However
$$\sum_{i, j=1}^{N} t_i t_j \frac{\left (|\langle u_i, u_j \rangle| - \varphi \|u_i\| \|u_j\| \right )^2}{\|u_i\| \|u_j\|}  =   \sum_{i, j = 1}^{N} t_i t_j \frac{|\langle u_i, u_j \rangle|^2}{\|u_i\| \|u_j\|} - 2 \varphi \sum_{i, j=1}^{N} t_i t_j |\langle u_i, u_j \rangle| + \varphi^2 \left ( \sum_{i=1}^{N} t_i \|u_i\| \right )^2.$$
Therefore, the estimate \eqref{mainestimate} rewrites as
$$\sum_{i, j=1}^{N} t_i t_j |\langle u_i, u_j \rangle| \leq \frac{1}{2\varphi} \sum_{i, j = 1}^{N} t_i t_j \frac{|\langle u_i, u_j \rangle|^2}{\|u_i\| \|u_j\|} + \frac{1}{2 \varphi} \left ( \varphi^2 - \frac{(1-\varphi)^2}{d} \right) \left ( \sum_{i=1}^{N} t_i \|u_i\| \right )^2.$$
The Cauchy--Schwarz inequality together with \eqref{eqnorma} yields now
\begin{align*}
\sum_{i, j=1}^{N} t_i t_j \frac{|\langle u_i, u_j \rangle|^2}{\|u_i\| \|u_j\|} &\leq \sqrt{\sum_{i, j=1}^{N} t^2_i \frac{|\langle u_i, u_j \rangle|^2}{\|u_i\|^2}} \cdot \sqrt{\sum_{i, j=1}^{N} t^2_j \frac{|\langle u_i, u_j \rangle|^2}{\|u_j\|^2}}\\
&=\sqrt{ \sum_{i=1}^{N} t_i^2} \cdot \sqrt{\sum_{j=1}^{N} t_j^2}=\|t\|^2=1.
\end{align*}
Moreover,
\begin{equation}
\label{cst}
\left( \sum_{i=1}^{N} t_i \|u_i\| \right)^2 \leq \|t\|^2 \left ( \sum_{i=1}^{N} \|u_i\|^2 \right ) = n
\end{equation}
by \eqref{eqtrace}. It is clear that equality holds in the last two estimates if and only if there exists $c>0$ such that $\|u_i\|=ct_i$ for every $1 \leq i \leq N$ (we note that $\langle u_i, u_j \rangle \neq 0$ by \eqref{condequi}). From \eqref{eqtrace} and the equality $\|t\|=1$ it follows that $c = \sqrt{n}$. Moreover, using the explicit formulas for $\varphi=\varphi_{\mathbb K}(n)$ and $d=d_{\mathbb K}(n)$, one easily checks that $\varphi^2 > \frac{(1-\varphi)^2}{d}$ for all $n \geq 2$. Because this coefficient is positive, applying the bounds yields
\begin{align*}
\sum_{i, j=1}^{N} t_i t_j |\langle u_i, u_j \rangle| &\leq \frac{1}{2\varphi} \sum_{i, j = 1}^{N} t_i t_j \frac{|\langle u_i, u_j \rangle|^2}{\|u_i\| \|u_j\|} + \frac{1}{2 \varphi} \left ( \varphi^2 - \frac{(1-\varphi)^2}{d} \right) \left ( \sum_{i=1}^{N} t_i \|u_i\| \right )^2\\
 &\leq \frac{1}{2 \varphi} + \frac{n}{2 \varphi} \left ( \varphi^2 - \frac{(1-\varphi)^2}{d} \right).
 \end{align*}
It is simple to check, by substituting the exact values for $\varphi=\varphi_{\mathbb{K}}(n)$ and for $d=d_{\mathbb{K}}(n)$, that the above expression is equal to $\delta_{\mathbb{K}}(n)$. Thus, the desired estimate follows. If $A_0$ is a non-empty set, then let $I = \{1, \ldots, N\} \setminus A_0$, with $|I|=N'$ be the set of indices for which vectors $u_i$ are non-zero and let $t' \in \mathbb{R}^{N'}$ be the corresponding vector defined as $t'_i=t_i$ for $i \in I$. Clearly, the vectors $u_i$ for $i \in I$ still form a tight frame in $\mathbb{K}^n$ with a constant $1$ and $\|t'\| \leq 1$. Therefore, the previous reasoning can still be applied for the vectors $u_i$ with $i \in I$ and the vector $t'$, yielding the same bound (we simply ignore the zero vectors $u_i$). However, if we had a strict inequality $\|t'\|<1$, then for example the estimate \eqref{cst} would be strict. In other words, if the equality holds in the main estimate, then we must have $\|t'\|=1$. This means that for $i \in A_0$, i.e. for $u_i=0$, we must also have $t_i=0$. This proves the first equality condition. It also establishes the fourth equality condition for those vectors $u_i$ that are equal to zero. For a non-zero $u_i$ the fourth condition follows immediately from the equality case in the Cauchy--Schwarz inequality that was used in the estimate \eqref{cst}. 

We are left with the other two equality conditions. The condition \eqref{condequi} shows that for every pair of vectors $\frac{u_i}{\|u_i\|}$ and $\frac{u_j}{\|u_j\|}$, either they span the same one-dimensional subspace of $\mathbb{K}^n$, or the modulus of their inner product is equal to $\varphi$. Let $w_1,\ldots,w_m\in\mathbb{K}^n$ be unit vectors obtained by choosing exactly one representative from each equivalence class of the relation 
$$ x\sim y \iff x=\alpha y$$
for some unimodular scalar $\alpha$. Then $w_1,\ldots,w_m$ form an equiangular set. For $1 \leq j \leq m$, let
$$A_j=\left \{\,i:\frac{u_i}{\|u_i\|}=\alpha w_j \text{ for some unimodular }\alpha \right \}.$$
 Let us also denote $S_j = \sum_{i \in A_j} t_i^2$ for $1 \leq j \leq m$. Since the vectors $u_i$ form a tight frame with constant $1$ we can now write
\begin{equation}
\label{etfeq2}
\begin{split}
x &= \sum_{i=1}^{N} \langle x, u_i \rangle u_i = \sum_{j=1}^{m} \left ( \sum_{i \in A_j}  \|u_i\|^2 \langle x, \alpha_i w_j \rangle (\alpha_i w_j) \right ) \\ &=
\sum_{j=1}^{m} \left ( \sum_{i \in A_j}  \|u_i\|^2 \right ) \langle x, w_j \rangle w_j =n \sum_{j=1}^{m} S_j \langle x, w_j \rangle w_j.
\end{split}
\end{equation}
In particular, $\id_n \in \lin \{ w_1w_1^*, \ldots, w_mw_m^* \}$. Moreover, since for any $i \in A_j$ we also have 
$$L_i = \sqrt{t_i} \|u_i\|^{\frac{-3}{2}} u_iu_i^{*} = \sqrt{t_i} \|u_i\|^{\frac{-3}{2}} \|u_i\|^2 (\alpha_i w_j) (\alpha_i w_j)^{*} = \sqrt{t_i} \|u_i\|^{\frac{1}{2}} w_j w_j^{*} $$
it follows that $\lin \{ w_1w_1^*, \ldots, w_mw_m^* \}=\mathcal{H}$. Our goal is to establish that $m = d= d_{\mathbb{K}}(n)$. Indeed, if $\dim \mathcal{H} < d$, then the estimate \eqref{trdim} becomes
$$\tr(A^2)=\tr(A^*A) \geq \frac{\tr(A)^2}{\rk A} \geq \frac{\tr(A)^2}{\dim \mathcal{H}} > \frac{\tr(A)^2}{d}.$$
This in turn leads to the fact that the main estimate is strict, contradicting the assumed equality. Therefore, we conclude that $m = d$. Since the vectors $w_1, \ldots, w_d$ form an equiangular set of the maximal possible cardinality, it follows that they necessarily form an ETF in $\mathbb{K}^n$ (see Theorem $5.10$ in \cite{foucart}). Thus, by comparing the traces, we deduce that they form a tight frame with constant $\frac{n}{d}$. In other words, for every $x \in \mathbb{K}^n$ we have
$$x = \frac{n}{d} \sum_{i=1}^{d} \langle x, w_i \rangle w_i.$$
Since we have $\dim \mathcal{H}=d$, the rank one matrices $w_1w_1^*, \ldots, w_dw_d^*$ are linearly independent and by comparing the equality above with \eqref{etfeq2}, it follows that $S_i=\frac{1}{d}$ for every $1 \leq i \leq m$. This proves the necessity of all conditions for the equality.

It remains to prove that the stated conditions are sufficient for equality. Assuming that all conditions are met, we calculate that for fixed $1 \leq k < l \leq d$ we have
\begin{align*}
\sum_{i \in A_k} \sum_{j \in A_l} t_it_j |\langle u_i, u_j \rangle| &= \sum_{i \in A_k} \sum_{j \in A_l} t_it_j \|u_i\| \|u_j\| |\langle w_k, w_l \rangle| = n \sum_{i \in A_k} \sum_{j \in A_l} t_i^2 t_j^2 \varphi \\
&=n \varphi \left ( \sum_{i \in A_k} t_i^2 \right ) \left ( \sum_{j \in A_l} t_j^2 \right ) = \frac{n \varphi}{d^2}.
\end{align*}
Similarly, for each $1 \leq k \leq d$ we obtain
\begin{align*}
\sum_{i \in A_k} \sum_{j \in A_k} t_it_j |\langle u_i, u_j \rangle| &= \sum_{i \in A_k} \sum_{j \in A_k} t_it_j \|u_i\| \|u_j\| |\langle w_k, w_k \rangle|\\
&=n\sum_{i \in A_k} \sum_{j \in A_k} t_i^2 t_j^2 =n \left ( \sum_{i \in A_k} t_i^2 \right )^2 = \frac{n}{d^2}.
\end{align*}
Therefore
\begin{align*}
\sum_{i, j=1}^{N} t_it_j |\langle u_i, u_j \rangle| = \sum_{k, l=1}^{d} \sum_{i \in A_k} \sum_{j \in A_l} t_it_j |\langle u_i, u_j \rangle| &=\frac{n \varphi}{d^2} \cdot d(d-1) + \frac{n}{d^2} \cdot d\\
&=\frac{n}{d} \left ( (d-1) \varphi +  1 \right ).
\end{align*}
A direct substitution of the formulas for $d=d_{\mathbb K}(n)$ and $\varphi=\varphi_{\mathbb K}(n)$ shows that the last expression is equal to $\delta_{\mathbb K}(n)$. Hence equality holds and the proof is complete.

\end{proof}

In the following remark we discuss notions of the maximal and quasimaximal relative projection constants and a consequence of Theorem \ref{twrbound}.

\begin{remark}
\label{remprojconst}
As noted already in the introductory section, for any $n \geq 1$ there exists an $n$-dimensional subspace $X$ of $\ell_{\infty}$ (over $\mathbb{K}$) such that $\lambda(X) = \lambda_{\mathbb{K}}(n)$. Therefore, for an integer $N \geq n$ it is natural to define a \emph{maximal relative projection constant} $\lambda_{\mathbb{K}}(n, N)$ as $\lambda_{\mathbb{K}}(n, N) = \sup \lambda(Y, \ell_{\infty}^N)$, where a supremum (in fact, a maximum) is taken over all $n$-dimensional subspaces $Y \subseteq \ell_{\infty}^N$. Maximal relative projection constants satisfy in particular $\lambda_{\mathbb{K}}(n, N) \leq \lambda_{\mathbb{K}}(n, N+1)$ and the value $\lambda_{\mathbb{K}}(n, N)$ converges to $\lambda_{\mathbb{K}}(n)$ as $N \to \infty$. Moreover, the formula \eqref{formula} in this case is exactly the same as before, but just with $N$ fixed, i.e.
\begin{equation}
\label{formula2}
\lambda_{\mathbb{K}}(n, N)=  \sup \sum_{i, j=1}^{N} t_i t_j |\langle u_i, u_j \rangle|
\end{equation}
where the supremum (maximum) again ranges over all non-negative vectors $t \in \mathbb{R}^N$ with $\|t\|=1$ and vectors $u_1, \ldots, u_N \in \mathbb{K}^n$ forming a tight frame with constant $1$ in $\mathbb{K}^n$. It is natural to ask whether the value of $\lambda_{\mathbb K}(n,N)$ stabilizes at $\lambda_{\mathbb K}(n)$ for sufficiently large $N$. This is equivalent to asking whether there always exists an $n$-dimensional polyhedral normed space with the maximal projection constant. Basso claimed this for $\mathbb K=\mathbb R$ in Theorem~1.4 of \cite{basso}, but the proof ultimately turned out to be incorrect. The author later published a corrigendum \cite{bassoerratum}, in which the mistake was explained in detail. Thus, it remains open whether, for every $n$, there exists $N_0$ such that \[ \lambda_{\mathbb K}(n,N)=\lambda_{\mathbb K}(n) \] for all $N\ge N_0$. Equivalently, it is unknown whether the supremum in \eqref{formula} can always be attained for some specific positive integer $N$. At present, the only situations in which we know that this happens, are exactly those for which the supremum in the estimate \eqref{formula} can be determined explicitly and it turns out that in those situations, it can be indeed attained for a specific $N$. It should be noted that Theorem \ref{twrchar} shows that, apart from the $\mathbb{R}^2$ case, there always exist non-polyhedral maximizers as well.

The upper bound $\lambda_{\mathbb{K}}(n)\leq \delta_{\mathbb{K}}(n)$ implies in particular that $\lambda_{\mathbb{K}}(n,N)\leq \delta_{\mathbb{K}}(n)$, and by Theorem \ref{twrbound} equality holds if and only if a maximal ETF exists in $\mathbb{K}^n$ and $N\ge d_{\mathbb{K}}(n)$. However, the maximal relative projection constants $\lambda_{\mathbb{K}}(n,N)$ can also be difficult to compute in general. For example, it was shown in \cite{chalmerslewicki} that $\lambda_{\mathbb{R}}(3, 5) = \frac{5+4\sqrt{2}}{7}$ but obtaining this value required substantial effort. One can reduce the difficulty of determining the supremum in \eqref{formula2} by restricting to the case where all weights $t_i$ are equal to $\frac1{\sqrt N}$. The resulting expression is called the \emph{quasimaximal relative projection constant} $\mu_{\mathbb{K}}(n,N)$, namely
$$\mu_{\mathbb{K}}(n, N) =\sup \frac{1}{N} \sum_{i, j=1}^{N} |\langle u_i, u_j \rangle|,$$
where the supremum is taken over all tight frames with constant $1$ in $\mathbb{K}^n$. Obviously $\mu_{\mathbb{K}}(n, N)$ gives a lower bound to $\lambda_{\mathbb{K}}(n, N)$ and is generally easier to handle, since one only needs to optimize over the choice of a tight frame. Deręgowska and Lewandowska, in their proof from \cite{deregowskalewandowska}, first established the inequality  $\mu_{\mathbb{K}}(n,N)\leq \delta_{\mathbb{K}}(n)$ for every $N$, and then showed that $\lambda_{\mathbb{K}}(n)$ can be approximated arbitrarily well by $\mu_{\mathbb{K}}(n,N)$ for suitable choices of $N$. For a general overview of the relations between maximal and quasimaximal relative projection constants, the reader is referred to \cite{foucartskrzypek}. 

We note that an immediate consequence of Theorem~\ref{twrbound} is the following observation: if a maximal ETF exists in $\mathbb{K}^n$, then the equality $\mu_{\mathbb{K}}(n,N)=\lambda_{\mathbb{K}}(n)$ holds if and only if $N$ is divisible by $d_{\mathbb{K}}(n)$. Indeed, if all coordinates of the vector $t$ are equal to $\frac1{\sqrt N}$, then condition $(3)$ becomes $\frac{|A_j|}{N} = \frac1{d_{\mathbb{K}}(n)}$ for every $1\le j\le d_{\mathbb{K}}(n)$. Hence $|A_j|=\frac{N}{d_{\mathbb{K}}(n)}$, so $N$ must be divisible by $d_{\mathbb{K}}(n)$. Conversely, if $N=kd_{\mathbb{K}}(n)$, then one can take $k$ copies of a fixed maximal ETF in $\mathbb{K}^n$ and verify directly that all equality conditions are satisfied.

\end{remark}

\section{Characterization of spaces with the maximal projection constant in the polyhedral case}
\label{seccharpolyhedral}

In this section, we establish the polyhedral case of Theorem~\ref{twrchar}. In other words, whenever a maximal ETF exists in $\mathbb{K}^n$, we characterize the spaces with the maximal projection constant within the class of polyhedral spaces. We begin with a simple lemma concerning an important property of maximal ETFs.

\begin{lem}
\label{lemwj}
Let $n\ge2$ be an integer such that there exists a maximal ETF $w_1,\ldots,w_d$ in $\mathbb K^n$ (of cardinality $d=d_{\mathbb K}(n)$). Then for every $j\in\{1,\ldots,d\}$ we have
$$w_j = C \sum_{i=1}^{d} \sgn \langle w_i, w_j \rangle w_i,$$
where $C=\frac{n}{d\delta_{\mathbb{K}}(n)}.$
\end{lem}

\begin{proof}
Let us fix $j \in \{1, \ldots, d\}$. From the fact that vectors $w_1, \ldots, w_d$ form an ETF in $\mathbb{K}^n$ it follows that
$$w_j = \frac{n}{d} \sum_{i=1}^{d} \langle w_j, w_i \rangle w_i = \frac{n}{d} \left ( w_j + \varphi_{\mathbb{K}}(n) \sum_{i \neq j} \sgn \langle w_i, w_j \rangle w_i \right ),$$
so that
$$\frac{d-n}{n\varphi_{\mathbb{K}}(n) }w_j =  \sum_{i \neq j} \sgn \langle w_i, w_j \rangle w_i.$$
Hence
$$C \sum_{i=1}^{d} \sgn \langle w_i, w_j \rangle w_i = Cw_j + \frac{C(d-n)}{n\varphi_{\mathbb{K}}(n)} w_j = \frac{d - n + n \varphi_{\mathbb{K}}(n) }{ d \delta_{\mathbb{K}}(n) \varphi_{\mathbb{K}}(n) } w_j.$$
Substituting the explicit formulas for $d=d_{\mathbb{K}}(n)$, $\delta_{\mathbb{K}}(n)$ and $\varphi_{\mathbb{K}}(n)$, it follows from a direct calculation that
$$\frac{d - n + n \varphi_{\mathbb{K}}(n) }{ d \delta_{\mathbb{K}}(n) \varphi_{\mathbb{K}}(n) }  = 1$$
and the proof is finished.
\end{proof}

Theorem \ref{twrchar} states that, after a suitable linear transformation, the dual unit ball of every space with the maximal projection constant is contained between the absolutely convex hull of a maximal ETF and an appropriately rescaled zonotope generated by the same vectors. In the following lemma we prove that the inclusion between these sets is non-strict only in the $\mathbb{R}^2$ case.

\begin{lem}
\label{lemincl}
Let $n \geq 2$ be an integer such that there exists a maximal ETF $w_1, \ldots, w_d$ in $\mathbb{K}^n$ (of cardinality $d = d_{\mathbb{K}}(n)$). Then
$$\absconv \{w_1, w_2, \ldots, w_d\} \subseteq \frac{n}{d\delta_{\mathbb{K}}(n)} Z(w_1, \ldots, w_d)$$
and the equality holds if and only if $n=2$ and $\mathbb{K}=\mathbb{R}$.
\end{lem}

\begin{proof} Let us denote $C=\frac{n}{d\delta_{\mathbb{K}}(n)}$, $K = \absconv \{w_1, w_2, \ldots, w_d\}$ and also $L=\frac{n}{d\delta_{\mathbb{K}}(n)} Z(w_1, \ldots, w_d)$. The inclusion $K \subseteq L$ follows directly from Lemma \ref{lemwj} for arbitrary dimension $n$ and field $\mathbb{K}$, as $w_j \in L$ for every $1 \leq j \leq d$. Hence also $K = \absconv \{w_1, \ldots, w_d \} \subseteq L$.

Now let us prove that in the $\mathbb{R}^2$ case we have the equality $K=L$. Here, $d=3$, $\delta_{\mathbb{R}}(2)=\frac{4}{3}$ and $w_1, w_2, w_3$ are vertices of some equilateral triangle inscribed in the unit circle. For $x \in L$
$$x =  a_1w_1 + a_2w_2 + a_3w_3,$$
where $|a_i| \in \left [ 0, \frac{1}{2} \right ]$ for $1 \leq i \leq 3$. Without loss of generality we can assume that $a_1 \geq a_2 \geq a_3$. From $w_1+w_2+w_3=0$ we have
$$x = (a_1-a_2)w_1 + (a_3-a_2)w_2,$$
and moreover
$$|a_1 - a_2| + |a_3-a_2|=(a_1-a_2)+(a_2-a_3)=a_1-a_3 \leq \frac{1}{2} + \frac{1}{2}=1,$$
so $x \in \absconv \{w_1, w_2, w_3 \}=K$. This proves that $L \subseteq K$, so $K=L$.

We are left with proving that in all other situations the inclusion is strict. Let $X=(\mathbb{K}^n, \| \cdot \|_X)$ and $Y=(\mathbb{K}^n, \| \cdot \|_Y)$ be normed spaces such that $B_{X^*}=K$ and $B_{Y^*}=L$. Hence, for $x \in \mathbb{K}^n$
$$\|x\|_X = \max_{1 \leq i \leq d} |\langle x, w_i \rangle|$$
and
$$\|x\|_Y=C \max_{|a_i| \leq 1} |\langle x, a_1w_1 + \ldots + a_dw_d \rangle | = C \sum_{i=1}^{d} | \langle x, w_i \rangle |.$$
The equality $B_{Y^*}=L = K=B_{X^*}$ is by duality equivalent to $B_X = B_Y$ and, in terms of the norms, it is equivalent to the fact that for every $x \in \mathbb{K}^n$ we have $\|x\|_Y = \|x\|_X$, i.e.
\begin{equation}
\label{eqnorms}
\max_{1 \leq i \leq d} |\langle x, w_i \rangle| =  C \sum_{i=1}^{d} | \langle x, w_i \rangle |.
\end{equation}
From the established inclusion $K \subseteq L$ it follows that $\|x\|_Y \geq \|x\|_X$ for every $x \in \mathbb{K}^n$ and our goal is to prove that this inequality is strict for some $x$.

We need to consider separately the case $n=2$ also for $\mathbb{K}=\mathbb{C}$. In this case, a maximal ETF (of cardinality $d=4$) is known to be unique up to unitary transformation of $\mathbb{C}^2$ (see for instance Example $3.2$ in \cite{hughstonsalamon}) and can be given as
$$\left \{ \left (1, 0 \right ), \frac{1}{\sqrt{3}}\left ( 1, \sqrt{2} \right ), \frac{1}{\sqrt{3}} \left (1, \sqrt{2}\omega \right ), \frac{1}{\sqrt{3}} \left  ( 1, \sqrt{2}\omega^2 \right )  \right\},$$
where $\omega \in \mathbb{C}$ satisfies $\omega^3=1$ and $\omega \neq 1$. Because $\delta_{\mathbb{C}}(2)=\frac{1+\sqrt{3}}{2}$ the inequality $\|x\|_Y > \|x\|_X$ takes the form
$$\sqrt{3}|x_1| + |x_1 + \sqrt{2}x_2| + |x_1 + \sqrt{2}\omega x_2| + |x_1 + \sqrt{2}\omega^2x_2|$$
$$> (1+\sqrt{3}) \max \left \{ \sqrt{3}|x_1|, |x_1 + \sqrt{2}x_2|, |x_1 + \sqrt{2}\omega x_2|, |x_1 + \sqrt{2}\omega^2x_2| \right \}$$
for $x=(x_1, x_2) \in \mathbb{C}^2$. This holds for example for $x=(0, 1)$, as $3\sqrt{2} > (1+\sqrt{3})\sqrt{2}.$

Assume now that $n>2$. For the sake of contradiction, let us suppose that for every $x \in \mathbb{K}^n$ we have $\|x\|_Y = \|x\|_X$. Setting $x=w_1$ gives
$$1=\|w_1\|_X = C \sum_{i=1}^{d} |\langle w_1, w_i\rangle|$$
(by Lemma \ref{lemwj} this actually holds for any $1 \leq i \leq d$, regardless of the assumed equality of the norms $\| \cdot \|_X, \| \cdot \|_Y$). Now, we shall consider real and complex cases separately, where the argument will be slightly different.

In the case $\mathbb{K}=\mathbb{R}$ let $v \in \mathbb{R}^n$ be a unit vector orthogonal to $w_1$. We can choose it in such a way that $v$ is not perpendicular to any vector $w_i$ for $i>1$. Indeed, vector $v$ is chosen from a subspace of dimension $n-1>1$ and hence we can choose it to be not orthogonal to any given finite set of non-zero vectors. Let us assume that index $k \in \{2, \ldots, d\}$ is such that 
$$|\langle v, w_k \rangle | = \max_{1 \leq i \leq d} | \langle v, w_i \rangle |=M.$$
By the Cauchy--Schwarz inequality we clearly have $M<1$ (as no vector $w_i$ is perpendicular to $w_1$). By multiplying $w_2, \ldots, w_d$ by appropriate unimodular scalars we can suppose that $\langle w_1, w_i \rangle = |\langle w_1, w_i \rangle | = \varphi$ for any $i \geq 2$ (where $\varphi=\varphi_{\mathbb{R}}(n)$). Similarly, by adjusting $v$ we can suppose also that $\langle v, w_k \rangle =|\langle v, w_k \rangle |=M$. We observe that if a real $t$ satisfies $0 \leq t \leq \frac{1-\varphi}{M}$, then $\|w_1-tv\|_X=1$, since for $2 \leq i \leq d$ we have
$$|\langle w_i, w_1 - tv \rangle | = | \varphi - t\langle w_i, v \rangle| \leq \varphi + tM \leq 1.$$
Now let us define a function $f: \left [0, \frac{1-\varphi}{M} \right ] \to \mathbb{R}$ as
$$f(t) = \sum_{i=2}^{d} |\langle w_i, w_1 - tv \rangle | = \sum_{i=2}^{d} |\langle w_i, v \rangle | \left |\frac{\langle w_i, w_1 \rangle}{\langle w_i, v \rangle } - t \right |$$
(we recall that $\langle w_i, v \rangle \neq 0$ for $i \geq 2$). From the assumed equality  \eqref{eqnorms} it now follows that for any $t \in \left [0, \frac{1-\varphi}{M} \right ]$ we have
\begin{equation}
\label{maxfunc}
\begin{split}
1 &=  \max_{1 \leq i \leq d} |\langle w_i, w_1-tv \rangle| = |\langle w_1, w_1-tv \rangle| \\
&= C  \sum_{i=1}^{d} | \langle w_i, w_1-tv \rangle | = C + Cf(t),
\end{split}
\end{equation}
so the function $f$ is a constant function equal to $\frac{1-C}{C}$. In particular, this real function must be differentiable at the point $\frac{\langle w_k, w_1 \rangle}{\langle w_k, v \rangle } = \frac{\varphi}{M} \in \left (0, \frac{1-\varphi}{M} \right )$. We note that this point is in the interior of the domain of $f$, as it is simple to verify that $\varphi_{\mathbb{R}}(n) < \frac{1}{2}$ for $n \geq 3$. If we now denote by $S \subseteq \{2, 3, \ldots, d\}$ the set of all indices $i$, for which the equality $\frac{\langle w_i, w_1 \rangle}{\langle w_i, v \rangle } = \frac{\varphi}{M}$ holds ($S$ is non-empty since $k \in S$), then the function $f$ can be written as
\begin{equation}
\label{derivative}
f(t) = A\left | \frac{\varphi}{M} - t \right | + \sum_{i \not \in S} |\langle w_i, v \rangle | \left |\frac{\langle w_i, w_1 \rangle}{\langle w_i, v \rangle } - t \right |,
\end{equation}
where $A = \sum_{i \in S} |\langle w_i, v \rangle | > 0 $. However, for a fixed $a \in \mathbb{R}$ the function $\mathbb{R} \ni t \to |a-t| \in \mathbb{R}$ is not differentiable at the point $a$. This implies that the right-hand side of the equality above is not differentiable at the point $\frac{\varphi}{M} \in \left (0, \frac{1-\varphi}{M} \right )$.  We obtained a contradiction, which proves that there exists a vector $x \in \mathbb{R}^n$ such that $\|x \|_Y > \|x \|_X$.

We are left with the complex case. Since the ETF vectors $w_1, \ldots, w_d$ span $\mathbb{C}^n$ and $n \geq 3$, without loss of generality we can assume that the vectors $w_1, w_2, w_3$ are linearly independent. By multiplying $w_2, \ldots, w_d$ by appropriate unimodular scalars we can assume that $\langle w_1, w_i \rangle = |\langle w_1, w_i \rangle | = \varphi$ for any $i \geq 2$ (where $\varphi=\varphi_{\mathbb{C}}(n)$). Since $w_1, w_2, w_3$ are linearly independent, there exists a vector $v \in \mathbb{C}^n$ satisfying the following system of linear equations:
$$\langle v, w_1 \rangle=0, \qquad \langle v, w_2 \rangle = 1, \qquad \langle v, w_3 \rangle = i.$$
Since $v \neq 0$, we can normalize $v$ into the unit vector. Moreover, by multiplying $v$ by a suitable unimodular scalar we can assume that similarly as in the previous case we have 
$$\langle v, w_k \rangle = \max_{1 \leq i \leq d} | \langle v, w_i \rangle | = M < 1$$ for some $k \in \{2, \ldots, d\}$. We note that after this modification of $v$, we still have $\langle v, w_2 \rangle = \lambda$ and $\langle v, w_3 \rangle = i \lambda$ for some non-zero $\lambda \in \mathbb{C}$. Thus, at least one of the numbers $\langle v, w_2 \rangle$ or $\langle v, w_3 \rangle$ remains non-real.

Calculating as in \eqref{maxfunc} we deduce that the function  $f: \left [0, \frac{1-\varphi}{M} \right ] \to \mathbb{R}$ defined as
$$f(t) = \sum_{i=2}^{d} |\langle w_i, w_1 - tv \rangle | = \sum_{i=2}^{d}|\langle w_i, w_1 \rangle - t \langle w_i, v \rangle |$$
is a constant function. This time we will examine the second derivative of $f$ on the smaller open interval $\left(0, \frac{\varphi}{M}\right)$. Note that since $\varphi_{\mathbb{C}}(n) \leq \frac{1}{2}$ for $n \ge 3$, this is indeed a sub-interval of $\left[0, \frac{1-\varphi}{M}\right]$. For every $i\ge2$ and $t \in \left(0, \frac{\varphi}{M}\right)$, we have $t|\langle w_i, v \rangle| \le tM < \varphi$. Thus, the complex number $\varphi - t\langle w_i, v \rangle$ is never zero, meaning the real function $t\mapsto |\varphi-t\langle w_i,v\rangle|$ is twice continuously differentiable on the entire interval $\left(0, \frac{\varphi}{M}\right)$. A straightforward computation now yields
$$ f''(t) = \sum_{i=2}^{d} \frac{\varphi^2\,\im(\langle w_i,v\rangle)^2} {\bigl(\varphi^2 -2\varphi\re(\langle w_i,v\rangle)t +|\langle w_i,v\rangle|^2t^2 \bigr)^{3/2}}.$$
Since $f$ is a constant function, $f''$ is identically equal to $0$ on this interval. Because $f''$ is a sum of non-negative terms, we must have $\im\langle w_i,v\rangle=0$ for every $2 \leq i \leq d$. However, by construction, at least one of the numbers $\langle v, w_2 \rangle$ or $\langle v, w_3 \rangle$ is non-real. This yields the desired contradiction and completes the proof.

\end{proof}

In the following lemma we establish a similar property like in Lemma \ref{lemwj}, but for vectors $u_i$ forming a tight frame with constant $1$ and attaining the equality in the estimate $\lambda_{\mathbb{K}}(n) \leq \delta_{\mathbb{K}}(n)$. 

\begin{lem}
\label{lemtiui}
Let $N \geq n \geq 2$ be integers and let $u_1, u_2, \ldots, u_N \in \mathbb{K}^n$ be vectors forming a tight frame with constant $1$. We assume that a vector $t \in \mathbb{R}^N$ with non-negative coordinates satisfies $\|t\|=1$ and
$$\sum_{i, j=1}^{N} t_it_j |\langle u_i, u_j \rangle| = \delta_{\mathbb{K}}(n).$$
Then for every $j \in \{1, \ldots, N\}$ we have
$$u_j =\frac{nt_j}{\delta_{\mathbb{K}}(n)} \sum_{i=1}^{N} t_i \sgn \langle u_i, u_j \rangle u_i.$$
\end{lem}

\emph{Proof.} Since the vectors $u_1, u_2, \ldots, u_N$ attain equality in Theorem \ref{twrbound}, all equality conditions have to be satisfied. We shall use the same notation as in statement of this theorem. Let us fix $j \in \{1, \ldots, N\}$. If $u_j=0$, then by the first condition of Theorem \ref{twrbound} we have also $t_j=0$ and hence there is nothing to prove. Otherwise, $u_j$ is a scalar multiple of one of the vectors of a maximal ETF $w_1, \ldots, w_{d}$ (where $d=d_{\mathbb{K}}(n)$) and $t_j \neq 0$. For each $k \in \{1, \ldots, d\}$ and $i \in A_k$ let us write $u_i = \sqrt{n}t_i \alpha_i w_k$ for some unimodular scalar $\alpha_i$. Without loss of generality let us assume that $j \in A_1$, that is $u_j=\sqrt{n} t_j \alpha_j w_1$. By combining the third condition of Theorem \ref{twrbound} with Lemma \ref{lemwj} we get
$$\sum_{i=1}^N t_i \sgn \langle u_i, u_j \rangle u_i =  \sum_{k=1}^{d} \sum_{i \in A_k} t_i \sgn \langle u_i, u_j \rangle u_i = \sum_{k=1}^{d} \sum_{i \in A_k} \sqrt{n} t^2_i \sgn \langle \alpha_i w_k, \alpha_j w_1 \rangle (\alpha_iw_k)$$
$$=\frac{\sqrt{n}}{d} \sum_{k=1}^{d} \overline{\alpha_i} \alpha_i \alpha_j \sgn \langle w_k, w_1 \rangle w_k = \frac{\sqrt{n}}{d} \cdot \frac{d \alpha_j \delta_{\mathbb{K}}(n)}{n}  w_1 = \frac{\alpha_j\delta_{\mathbb{K}}(n)}{\sqrt{n}} w_1 = \frac{\delta_{\mathbb{K}}(n)}{nt_j} u_j$$
and the conclusion follows. \qed

The following lemma lies at the core of the formula \eqref{formula} for the absolute projection constant. Various versions of this result are known in the literature, and an operator $E$ satisfying conditions of the type stated below is usually called a \emph{Chalmers--Metcalf operator}. The existence of such operators follows from trace duality. For a systematic treatment of Chalmers--Metcalf operators and their applications, we refer the reader to \cite{lewickiskrzypek}.

\begin{lem}
\label{lemcm}
Let $N \geq n \geq 1$ be integers. Assume that $X$ is an $n$-dimensional subspace of $\ell_{\infty}^N$ (considered over $\mathbb{K}$). Then there exists a linear operator $E: \ell_{\infty}^{N} \to \ell_{\infty}^N$ with $E(X) \subseteq X$ satisfying 
$$\tr (E|_X) = \lambda(X) \quad \text{ and } \quad \sum_{i=1}^{N} \|E(e_i)\|_{\infty}=1.$$
Furthermore, for any linear operator $E: \ell_{\infty}^{N} \to \ell_{\infty}^N$ satisfying $E(X) \subseteq X$ and $\sum_{i=1}^{N} \|E(e_i)\|_{\infty}=1$ we have $\left |\tr (E|_X) \right | \leq \lambda(X)$.
\end{lem}

\begin{proof}

Throughout the proof, we assume that $P: \ell_{\infty}^N \to X$ is a linear projection onto $X$ with $\|P\|=\lambda(X)$. 

We begin with the proof of the inequality. Let $E: \ell_{\infty}^N \to \ell_{\infty}^N$ be any linear operator satisfying $E(X) \subseteq X$ and $\sum_{i=1}^N \|E(e_i)\|_{\infty} = 1$. Since $E(X)\subseteq X$, the composition $E\circ P$ satisfies $(E\circ P)(\ell_{\infty}^N)\subseteq X $ and $ (E\circ P)|_X=E|_X$. Since the trace of a finite-rank operator coincides with the trace of its restriction to any invariant subspace containing its range, it follows that $\tr(E\circ P)=\tr(E|_X)$. Thus, by the cyclic property of the trace we also have $\tr(P \circ E) = \tr(E|_X)$. By definition, the trace of $P \circ E$ on $\ell_{\infty}^N$ can be calculated as:
$$\tr(P \circ E) = \sum_{i=1}^N (P(E(e_i)))_{i}.$$
Therefore
\begin{equation}
\label{tracecm}
\begin{split}
\left | \tr(E|_X) \right | &= \left | \tr(P \circ E) \right  | \leq \sum_{i=1}^N \left |(P(E(e_i)))_{i}\right | \leq \sum_{i=1}^{N} \left \| P(E(e_i)) \right \|_{\infty} \\
&\leq \sum_{i=1}^{N} \|P\| \|E(e_i)\|_{\infty} = \|P\| \sum_{i=1}^{N} \|E(e_i)\|_{\infty} = \lambda(X),
\end{split}
\end{equation}
which proves the inequality part of the lemma.

For the second part, we use the existence of a Chalmers--Metcalf operator for a subspace $X$. By Theorem~2.15 of \cite{lewickiskrzypek} there exists an operator $E: X \to \ell_{\infty}^N$, which satisfies $E(X) \subseteq X$ and 
$$\tr(E|_X) = \tr(E \circ P) = \|P\| = \lambda(X).$$ Furthermore, $E$ can be represented as a convex combination of rank-one operators formed by extremal pairs for $P$, that is 
$$E = \sum_{k=1}^m \alpha_k x_k \otimes f_k,$$
where 
\begin{itemize}
\item $\alpha_k > 0$ for $k=1, \ldots, m$ and $\sum_{k=1}^m \alpha_k = 1$, 
\item $x_k \in \ell_{\infty}^N$ and $\|x_k\|_{\infty}=1$ for $k=1, \ldots, m$,
\item $f_k \in X^*$ and $\|f_k\|=1$ for $k=1, \ldots, m$,
\item $f_k(P(x_k))=\|P\|$ for any $k=1, \ldots, m$.
\end{itemize}
To define an operator on the entire space $\ell_{\infty}^N$, we apply the Hahn--Banach theorem to extend each functional $f_k \in X^*$ to a functional $\tilde{f}_k \in (\ell_{\infty}^N)^*$ without increasing its norm, so that $\left \|\tilde{f}_k \right \|= 1$. Since the space $(\ell_{\infty}^N)^*$ is isometric to $\ell_1^N$, we can represent each extended functional as a sequence $\tilde{f}_k = (f_{k,1}, \ldots, f_{k,N})$ satisfying $\left \| \tilde{f}_k \right \|_1 = \sum_{i=1}^N |f_{k,i}| = 1$.

We define the linear operator $\widetilde{E}: \ell_{\infty}^N \to \ell_{\infty}^N$ by $\widetilde{E} = \sum_{k=1}^m \alpha_k x_k \otimes \tilde{f}_k$. For any $x \in X$ we have 
$$\widetilde{E}(x) = \sum_{k=1}^m \alpha_k x_k \tilde{f}_k(x) = \sum_{k=1}^m \alpha_k x_k f_k(x) = E(x).$$ 
Therefore, the restriction of $\widetilde{E}$ to $X$ is equal to $E$, which gives us $\widetilde{E}(X) \subseteq X$ and $\tr(\widetilde{E}|_X) = \tr(E|_X) = \lambda(X)$. Now, for any $i \in \{1, \ldots, N\}$ it follows that
$$\widetilde{E}(e_i) = \sum_{k=1}^m \alpha_k \tilde{f}_k(e_i) x_k = \sum_{k=1}^m \alpha_k f_{k,i} x_k,$$
and hence
$$\|\widetilde{E}(e_i)\|_{\infty} \leq \sum_{k=1}^{m} \alpha_k |f_{k, i}| \| x_k \|_{\infty} = \sum_{k=1}^{m} \alpha_k |f_{k, i}|.$$
Summing this over all $i = 1, \dots, N$ we obtain
$$\sum_{i=1}^N \|\widetilde{E}(e_i)\|_{\infty} \leq \sum_{i=1}^N \sum_{k=1}^m \alpha_k |f_{k,i}| = \sum_{k=1}^m \alpha_k \sum_{i=1}^N |f_{k,i}| = \sum_{k=1}^m \alpha_k = 1.$$
This proves that $\sum_{i=1}^N \|\widetilde{E}(e_i)\|_{\infty} \leq 1$. To show that we actually have an equality, it remains to notice that if we had $\sum_{i=1}^N \|\widetilde{E}(e_i)\|_{\infty} < 1$, then by the inequality chain obtained in \eqref{tracecm} we would get
$$|\tr(E|_X)|=|\tr(\widetilde{E}|_X)| < \lambda(X),$$
which contradicts the assumption about the operator $E$. This concludes the proof.

\end{proof}

K\"onig and Tomczak-Jaegermann claimed in Lemma~2.1 of \cite{konigtomczak} that a result analogous to the previous lemma holds for finite-dimensional subspaces of the infinite-dimensional $\ell_{\infty}$ space. However, in the following remark we present an example showing that such a result is false in general. The example was suggested by an anonymous referee.

\begin{remark}
\label{remarkcounterexample}
Let $v = (1, 1, \ldots) \in \ell_{\infty}$ and let $X = \lin \{v \}$ be a $1$-dimensional subspace. Then clearly $\lambda(X)=1$. Let $f \in (\ell_{\infty})^*$ be any Banach limit, that is, a norm-one functional satisfying satisfying $f(v)=1$ and $f(e_i)=0$ for $i \geq 1$. Consider an operator $E: \ell_{\infty} \to \ell_{\infty}$ defined for $x \in \ell_{\infty}$ by 
$$E(x) = (x_1+f(x))v.$$
Then $E(X) \subseteq X$, $E(e_1)=v$ and $E(e_i)=0$ for $i \geq 2$. Consequently, we have $\sum_{i=1}^{\infty} \|E(e_i)\|_{\infty}=1$. On the other hand, $E|_X=2 \id_X$, and therefore
$$\tr(E|_X)=2 > 1 = \lambda(X).$$
\end{remark}

In their proof of \eqref{formula}, K\"onig and Tomczak-Jaegermann introduced an auxiliary inner product $\langle \cdot, \cdot \rangle_E$ on a subspace $X$ of $\ell_{\infty}^N$, defined with the help of the operator $E$ as in the previous lemma. However, a priori it is only clear that this auxiliary inner product is merely a semi-inner product, i.e. it could happen that $\langle x, x \rangle_E=0$ for some non-zero vector $x \in X$. This difficulty can be overcome using the fact that the function $\delta_{\mathbb K}(n)$ is increasing in $n$. To address this issue, we shall use the following elementary algebraic lemma, whose proof is omitted.

\begin{lem}
\label{lemtrace}
Let $X$ be a finite-dimensional vector space over $\mathbb{K}$ and let $V \subseteq X$ be a linear subspace. Suppose that $E: X \to X$ is a linear operator satisfying $E|_V \equiv 0$. Let $\widetilde{E}:X/V \to X/V$ be the induced operator on the quotient space $X/V$. Then $\tr E = \tr \widetilde{E}$.
\end{lem}

The proof of the polyhedral case of Theorem~\ref{twrchar} combines ideas from both papers \cite{konigtomczak} and \cite{konigtomczak2} of K\"onig and Tomczak-Jaegermann.

\begin{proof}[Proof of Theorem \ref{twrchar} in the polyhedral case] Throughout the proof of both implications we assume that $X$ is a polyhedral space. We first note that the inclusions from the statement are equivalent to the fact that $X$ is isometric to a space whose norm is given by
\begin{equation}
\label{normateza}
\|x\| = \max \left \{ \max_{1 \leq i \leq d} | \langle x, w_i \rangle |, \: \max_{w \in I} |\langle x, w \rangle| \right \},
\end{equation}
for some finite subset (possibly empty) $I$ of the zonotope $\frac{n}{d\delta_{\mathbb{K}}(n)} Z(w_1, \ldots, w_d)$.

Let us first suppose that $X$ is an $n$-dimensional polyhedral space satisfying $\lambda(X) = \lambda_{\mathbb{K}}(n)$. In this case, we shall prove that $X$ is isometric to a space with the norm as above. Since there exists a maximal ETF in $\mathbb{K}^n$, Theorem \ref{twrbound} implies that $\lambda(X)=\lambda_{\mathbb{K}}(n)=\delta_{\mathbb{K}}(n)$. Since every $n$-dimensional polyhedral normed space has an isometric embedding into $\ell^N_{\infty}$ for some positive integer $N$ and the absolute projection constant is invariant under linear isometries, we may assume that $X$ is a subspace of $\ell_\infty^N$ satisfying $\lambda(X) = \delta_{\mathbb{K}}(n)$. Let $E:\ell_\infty^N\to\ell_\infty^N$ be an operator provided by Lemma~\ref{lemcm}, that is $E(X) \subseteq X$, $\sum_{i=1}^{N} \|E(e_i)\|_{\infty} = 1$ and $\tr(E|_X)=\delta_{\mathbb{K}}(n)$. For $1\le i\le N$, define $t_i = \sqrt{\|E(e_i)\|_{\infty}}$ so that
$$\sum_{i=1}^{N} t_i^2 = 1 \quad \text{ and } \quad |E(e_i)_j| \leq \|E(e_i)\|_{\infty} = t_i^2 \quad \text { for all } i, j \in \{1, 2, \ldots, N\}.$$

Let $S \subseteq \{1, 2, \ldots, N\}$ be the set of all indices $i \in \{1, \ldots, N\}$ such that $E(e_i) \neq 0$ (or equivalently $t_i \neq 0$) and let $W \subseteq \ell_{\infty}^N$ be a linear subspace defined as $W = \lin \{e_i: \ i \in \{1, \ldots, N\} \setminus S\}.$ Clearly, $E|_W \equiv 0$. Let us also take $V = W \cap X \subseteq X$ and let $X'=X/V$ be the quotient space. By Lemma~\ref{lemtrace} we have
\begin{equation}
\label{trace}
\tr(E|_X) = \tr(\widetilde{E}),
\end{equation}
where $\widetilde{E}: X' \to X'$ is the operator induced by $E|_X$ on the quotient space. In the space $X'$ we introduce an inner product $\langle \cdot, \cdot \rangle_E$ defined for $x, y \in X$ as
$$\langle x, y \rangle_E = \sum_{i=1}^{N} t_i^2 x_i \overline{y_i} .$$
If $x_1-x_2, y_1-y_2 \in V$ then $\langle x_1, y_1 \rangle_E = \langle x_2, y_2 \rangle_E$ and therefore $\langle \cdot, \cdot \rangle_E$ does not depend on the representation in $X'$. Moreover, we have $\langle x, x \rangle_E = 0$ if and only if $x \in V$. Hence, the function $\langle \cdot, \cdot \rangle_E$ is indeed a correctly defined inner product on $X'$.

Let $0 \leq m \leq n$ be the dimension of $X'$. If $m=0$ then by (\ref{trace}) we have 
$$0 = \tr(\widetilde{E}) = \tr(E|_X) = \lambda(X),$$
which is an obvious contradiction. Hence $m \geq 1$ and we may pick an orthonormal basis of $X'$ with respect to $\langle \cdot, \cdot \rangle_E$. Let this orthonormal basis be given by vectors $x_1, x_2, \ldots, x_m \in X$. Then again by (\ref{trace}) we have
$$\delta_{\mathbb{K}}(n)=\lambda(X)=\tr(E|_X) = \tr(\widetilde{E}) = \sum_{i=1}^{m} \langle E(x_i), x_i \rangle_E.$$
Now for any fixed $1 \leq i \leq m$ we can write
$$\langle E(x_i), x_i \rangle_E = \sum_{j=1}^{N} t_j^2 E(x_i)_j \overline{x_{ij}}  = \sum_{j=1}^{N} t_j^2 E \left ( \sum_{k=1}^{N} x_{ik}e_k \right )_j \overline{x_{ij}} =\sum_{j, k=1}^N t_j^2 E(e_k)_jx_{ik}\overline{x_{ij}}.$$
Hence
$$\sum_{i=1}^{m} \langle E(x_i), x_i \rangle_E = \sum_{i=1}^{m} \sum_{j, k=1}^N t_j^2 E(e_k)_jx_{ik}\overline{x_{ij}} = \sum_{j, k=1}^N \sum_{i=1}^{m} t_j^2 E(e_k)_jx_{ik}\overline{x_{ij}}.$$
For fixed $1 \leq j, k \leq N$ from the inequality $|E(e_k)_j| \leq t_k^2$ it now follows that
\begin{equation}
\label{oszacowanie}
\left | \sum_{i=1}^{m} t_j^2 E(e_k)_jx_{ik}\overline{x_{ij}} \right | = t_j^2 |E(e_k)_j| \left | \sum_{i=1}^{m} x_{ik}\overline{x_{ij}} \right | \leq t_j^2 t_k^2 \left | \sum_{i=1}^{m} x_{ik}\overline{x_{ij}} \right |.
\end{equation}
Summarizing we have
\begin{equation}
\label{oszstala}
\begin{aligned}
\delta_{\mathbb{K}}(n) & = \sum_{i=1}^{m} \langle E(x_i), x_i \rangle_E  \leq  \sum_{j, k=1}^N  t_j^2 t_k^2 \left | \sum_{i=1}^{m} x_{ik}\overline{x_{ij}} \right | & \\
& =  \sum_{j, k=1}^N t_j t_k \left | \sum_{i=1}^{m} t_kx_{ik}\overline{t_jx_{ij}} \right |. 
\end{aligned}
\end{equation}
Let us now introduce vectors $u_1, \ldots, u_N \in \mathbb{K}^m$ defined as $u_{ij}=t_i \overline{x_{ji}}$ for $1 \leq i \leq N$ and $1 \leq j \leq m$. Since the vectors $x_1, \ldots, x_m$ form an orthonormal basis with respect to $\langle \cdot, \cdot \rangle_E$, for any $1 \leq j \leq m$ we have
$$\sum_{i=1}^{N} |u_{ij}|^2 = \sum_{i=1}^{N} t_i^2 |x_{ji}|^2 = 1$$
and for $1 \leq j < k \leq m$
$$\sum_{i=1}^{N} u_{ij}\overline{u_{ik}} = \sum_{i=1}^{N} t_i^2 \overline{x_{ji}}x_{ki} = 0.$$
This shows that the vectors $u_1, \ldots, u_N$ form a tight frame in $\mathbb{K}^m$ with constant $1$. Thus, the conditions of Theorem \ref{twrbound} are satisfied and in consequence we have the estimate
$$\sum_{k, j=1}^N t_k t_j | \langle u_k,  u_j\rangle| \leq \delta_{\mathbb{K}}(m),$$
where $\langle \cdot, \cdot \rangle$ is a standard inner-product in $\mathbb{K}^m$. On the other hand, the estimate (\ref{oszstala}) reads as
$$\delta_{\mathbb{K}}(n) \leq \sum_{k, j=1}^N t_k t_j | \langle u_k,  u_j\rangle|,$$
and therefore $\delta_{\mathbb{K}}(n) \leq \delta_{\mathbb{K}}(m)$. However, since $m \leq n$ and it is simple to check that the function $\delta_{\mathbb{K}}(n)$ is strictly increasing in $n$, we must have $m=n$. It follows that $V = \{0\}$, i.e. the function $\langle \cdot, \cdot \rangle_E$ is in fact an inner product on the whole subspace $X$ and the vectors $x_1, \ldots, x_n$ form a linear basis of $X$. Let us note that up to this point we have only used $\lambda(X)>\delta_{\mathbb{K}}(n-1)$, which will be relevant for the proof in the general case. In particular, the norm of an arbitrary vector $\sum_{i=1}^{n} c_ix_i \in X$ (where $c_i \in \mathbb{K}$) is now given as
\begin{equation}
\label{norma}
\left \| \sum_{i=1}^{n} c_ix_i \right \|_{\infty} = \max_{1 \leq j \leq N} \left | \sum_{i=1}^{n} c_ix_{ij} \right |.
\end{equation}
Since the equality
$$\delta_{\mathbb{K}}(n) = \sum_{k, j=1}^N t_k t_j | \langle u_k,  u_j\rangle|,$$
holds, all equality conditions from Theorem \ref{twrbound} must be satisfied. Let a maximal ETF $w_1, w_2, \ldots, w_d \in \mathbb{K}^n$ and a partition $\{1, 2, \ldots, N\} = A_0 \sqcup A_1 \sqcup \ldots \sqcup A_d$ be as in the statement of Theorem \ref{twrbound} (in particular $S = \{1, 2, \ldots, N\} \setminus A_0 $). Besides the equality conditions of Theorem \ref{twrbound}, it is crucial to observe that the equality must also hold in the estimate (\ref{oszacowanie}). Therefore, for any fixed $j, k \in \{1, 2, \ldots, N\}$ we have
$$\left | \sum_{i=1}^{n} t_j^2 E(e_k)_jx_{ik}\overline{x_{ij}} \right | = t_j^2 |E(e_k)_j| \left | \sum_{i=1}^{n} x_{ik}\overline{x_{ij}} \right | = t_j^2 t_k^2 \left | \sum_{i=1}^{n} x_{ik}\overline{x_{ij}} \right |.$$
In particular, if $j, k \in S$ we have $t_j, t_k \neq 0$ and moreover
$$\left | \sum_{i=1}^{n} x_{ik}\overline{x_{ij}} \right | =\frac{1}{t_kt_j}\left | \sum_{i=1}^{n} (t_kx_{ik})\overline{(t_jx_{ij})} \right | = \frac{1}{t_kt_j}\left | \sum_{i=1}^{n} \overline{u_{ki}}u_{ji} \right | = \frac{| \langle u_j, u_k \rangle |}{t_kt_j} \neq 0,$$
because the equality conditions given in Theorem \ref{twrbound} imply that
 $$\frac{| \langle u_j, u_k \rangle |}{t_kt_j} = n \left | \left \langle  \frac{u_j}{\|u_j\|},  \frac{u_k}{\|u_k\|} \right \rangle \right | \in \{n, n \varphi_{\mathbb{K}}(n)\}.$$
Hence, since the equality holds in the estimates \eqref{oszacowanie} and \eqref{oszstala} it follows that for a fixed $k \in S$ we have $E(e_k)_j = \alpha_k \sgn \langle u_j, u_k \rangle t^2_k$ for all $j \in S$ and some scalar $\alpha_k$ of modulus $1$. Furthermore, by a similar calculation as before, we obtain
$$\delta_{\mathbb{K}}(n) = \sum_{i=1}^{n} \langle E(x_i), x_i \rangle_E = \sum_{j, k=1}^N E(e_k)_j \sum_{i=1}^{n} t_j^2  x_{ik}\overline{x_{ij}} = \sum_{j, k \in S} E(e_k)_j \sum_{i=1}^{n} t_j^2  x_{ik}\overline{x_{ij}}$$
$$= \sum_{j, k \in S} \alpha_k t_k t_j \sgn \langle u_j, u_k \rangle \langle u_j, u_k \rangle = \sum_{j, k \in S} \alpha_k t_k t_j |\langle u_j, u_k \rangle| = \sum_{j, k \in S} t_k t_j |\langle u_j, u_k \rangle|. $$
Since $|\alpha_k|=1$ this immediately implies that $\alpha_k=1$ for all $k \in S$. We can now summarize our knowledge about the operator $E$ as follows:
\begin{enumerate}
\item $E(e_k)_j = \sgn \langle u_j, u_k \rangle t^2_k$ for $j, k \in S$.
\item $E(e_k)_j = 0$ for $k \in \{1, 2, \ldots, N \} \setminus S$ and $j \in \{1, 2, \ldots, N\}$.
\item $|E(e_k)_j| \leq t_k^2$ for $k \in S$ and $j \in \{1, 2, \ldots, N\} \setminus S$.
\end{enumerate}
Thus, we do not know the exact values of $E(e_k)_j$ for $k \in S$, $j \in \{1, 2, \ldots, N\} \setminus S$ but we only have an estimate for their moduli. These entries can be non-zero in general and this is precisely the reason behind the possible appearance of additional points of the zonotope in $B_{X^*}$ beyond those coming from the absolutely convex hull of the vectors $w_1,\ldots,w_d$. For an illustration, let us consider some specific example. Suppose that $n=2$, $\mathbb{K}=\mathbb{R}$, $N=5$, $t_1=t_2=t_3 = \frac{1}{\sqrt{3}}$ and $t_4=t_5=0$. Then, taking into the account the sign pattern of three equiangular lines in $\mathbb{R}^2$, the matrix of the operator $E$ has the following form
$$\left [%
\begin{array}{ccccc}
  \frac{1}{3} & -\frac{1}{3} & -\frac{1}{3} & 0 & 0  \\
  -\frac{1}{3} & \frac{1}{3} & -\frac{1}{3} & 0 & 0 \\
  -\frac{1}{3} & -\frac{1}{3} & \frac{1}{3} & 0 & 0 \\
  a_1 & a_2 & a_3 & 0 & 0 \\
  b_1 & b_2 & b_3 & 0 & 0
\end{array}%
\right],
 $$
where $a_i, b_i \in \left [ -\frac{1}{3}, \frac{1}{3} \right ]$ are real numbers.

We shall now prove that these conditions on the operator $E$ imply that $E|_X = C \id_X$, where $C=\frac{\delta_{\mathbb{K}}(n)}{n}.$ Let us fix $1 \leq i \leq  n$. By the first two properties of the operator $E$ and Lemma \ref{lemtiui}, for any $j \in S$ we have
\begin{align*}
E(x_i)_j = \sum_{k=1}^{N} x_{ik} E(e_k)_j &= \sum_{k=1}^{N} x_{ik}t_k^2 \sgn \langle u_j, u_k \rangle=\sum_{k=1}^{N} t_k \overline{u_{ki} \sgn \langle u_k, u_j \rangle} \\
&= \frac{C}{t_j}\overline{u_{ji}} = Cx_{ij}.
\end{align*}
Hence $E(x_i)_j = C x_{ij}$ for any $j \in S$. Since $E(X) \subseteq X$ we can represent the vector $E(x_i) \in X$ in a basis $x_1, x_2, \ldots, x_n$ of $X$ as $\sum_{k=1}^{n} c_kx_k$ for some $c_k \in \mathbb{K}$. Then for $j \in S$ we have
$$\frac{C}{t_j}\overline{u_{ji}} = C x_{ij} =   E(x_i)_j =  \sum_{k=1}^{n} c_k x_{kj} = \sum_{k=1}^{n} \frac{c_k}{t_j} \overline{u_{jk}},$$
which gives us
$$C u_{ji} = \sum_{k=1}^{n} c_k u_{jk}.$$
In fact, this equality holds for all $1 \leq j \leq N$, as for $j \not \in S$ we have $u_j=0$. However, because the vectors $u_1, u_2, \ldots, u_N$ form a tight frame with constant $1$ in the space $\mathbb{K}^n$, the $n$ vectors of the form $(u_{1k}, u_{2k}, \ldots, u_{Nk}) \in \ell_2^N$ for $1 \leq k \leq n$ are linearly independent, as they are pairwise orthogonal and non-zero (since they are unit in the Euclidean norm). In particular, by the linear independence of these vectors, it follows that $c_i=C$ and $c_k = 0$ for $k \neq i$. Thus $E(x_i) = Cx_i$ and the equality $E(x_i)_j = C x_{ij}$ is true for all $1 \leq j \leq N$. For an index $j \not \in S$ this leads us to
$$Cx_{ij} =  E(x_i)_j =  \sum_{k=1}^{N} x_{ik} E(e_k)_j = \sum_{k \in S} x_{ik} E(e_k)_j = \sum_{k \in S} \frac{\overline{u_{ki}}}{t_k} E(e_k)_j$$
$$=\sum_{s=1}^{d} \sum_{k \in A_s} \frac{\overline{u_{ki}}}{t_k} E(e_k)_j = \sqrt{n} \sum_{s=1}^{d} \sum_{k \in A_s} \overline{w_{si}} E(e_k)_j = \sqrt{n} \sum_{s=1}^d \overline{\gamma}_{s}\overline{w_{si}},$$
where $\overline{\gamma}_s = \sum_{k \in A_s} \overline{\alpha_k} E(e_k)_j$ for $1 \leq s \leq d$ and $\alpha_k$ is a unimodular scalar such that $u_k = \sqrt{n} t_k \alpha_k w_s$ for $k \in A_s$. We note that from properties $(2)$ and $(3)$ of operator $E$ we can estimate
$$|\gamma_s| \leq \sum_{k \in A_s} |E(e_k)_j| \leq \sum_{k \in A_s} t_k^2 = \frac{1}{d},$$
where the last equality follows from the equality condition $(3)$ in Theorem \ref{twrbound}. Thus, for $j \not \in S$ and an arbitrary vector $x = c_1x_1 + \ldots + c_nx_n \in X$ we have
$$\left | \sum_{i=1}^{n} c_ix_{ij} \right | = \frac{\sqrt{n}}{C} \left | \sum_{i=1}^{n} c_i \left ( \sum_{s=1}^d \overline{\gamma}_{s}\overline{w_{si}} \right ) \right | = \sqrt{n} \left | \langle c, w \rangle \right |,$$
where $c=(c_1, \ldots, c_n) \in \mathbb{K}^n$ and $w = \frac{1}{dC}(a_1w_1 + \ldots + a_dw_d) \in \frac{n}{d\delta_{\mathbb{K}}(n)} Z(w_1, \ldots, w_d)$ with $a_i = d \gamma_i$ (then $|a_i| \in [0, 1]$). Therefore, comparing this with \eqref{norma} we see that the indices $j \not \in S$ contribute to the part of the maximum in the desired formula \eqref{normateza} for the norm corresponding to additional points of the zonotope (up to rescaling by $\sqrt{n}$). On the other hand, if $j \in S$ then $u_j$ is a multiple of one of the vectors $w_1, \ldots, w_d$. Without loss of generality, let us assume that $u_j=\sqrt{n}t_j \alpha_j w_1$ for some unimodular $\alpha_j$. Then
$$\left | \sum_{i=1}^{n} c_ix_{ij} \right | = \left | \sum_{i=1}^{n} \frac{c_i\overline{u_{ji}}}{t_j} \right | = \sqrt{n} \left | \sum_{i=1}^{n} c_i\alpha_j\overline{w_{1i}} \right | = \sqrt{n} | \langle c, w_1 \rangle |.$$
Again, after rescaling by $\sqrt{n}$, this is a part of the desired formula \eqref{normateza}. Clearly, for each $1 \leq i \leq d$ there exists at least one vector $u_j$ which is a multiple of $w_i$ (third condition in Theorem \ref{twrbound}), so indeed, for every $1 \leq i \leq d$ an inner product $|\langle c, w_i \rangle|$ will appear in the formula for the norm. As explained before, the indices $\{1, \ldots, N \} \setminus S$ correspond exactly to additional points of the zonotope \eqref{normateza}. This completes the proof of the first implication.

For the converse implication, let us suppose that the norm of $X$ is given as in \eqref{normateza}. Since $X$ is assumed to be polyhedral, the set $I$ in the formula for the norm is finite. Let $|I|=m$, where $m \geq 0$ and $N=d+m$ (then $d \leq N$). Thus, we have
$$\|x\| =  \max \left \{ \max_{1 \leq i \leq d} | \langle x, w_i \rangle |, \: \max_{1 \leq i \leq m} |\langle x, v_i \rangle| \right \},$$
where $v_{ij} = \sum_{k=1}^{d} a_{ik}w_{kj}$ for $1 \leq i \leq m$, $1 \leq j \leq n$ and some scalars $a_{ik} \in \mathbb{K}$ satisfying $|a_{ik}| \leq \frac{n}{d\delta_{\mathbb{K}}(n)}$ (if $m=0$ then simply there is no second part in the formula). We consider an isometric embedding $\widetilde{X}$ of $X$ into $\ell_{\infty}^N$ given as
$$\mathbb{K}^n \ni x \to (\langle x, w_1 \rangle, \ldots \langle x, w_d \rangle, \langle x, v_1 \rangle, \ldots, \langle x, v_m \rangle ) \in \ell_{\infty}^N.$$
A linear basis $x_1, \ldots, x_n$ of $\widetilde{X}$ is thus given as $x_{ij} = \overline{w_{ji}}$ if $1 \leq j \leq d$ and $x_{i(d+j)} = \sum_{k=1}^{d} \overline{a_{jk}w_{ki}}$ for $1 \leq j \leq m$ (where $1 \leq i \leq n$).

Since $\lambda(X) \leq \delta_{\mathbb{K}}(n)$, to prove that $\lambda(X)=\lambda(\widetilde{X}) = \delta_{\mathbb{K}}(n)$, it is enough by Lemma \ref{lemcm} to construct an operator $E: \ell_{\infty}^N \to \ell_{\infty}^N$ such that $E(\widetilde{X}) \subseteq \widetilde{X}$, $\sum_{i=1}^{N} \|E(e_i)\|_{\infty}=1$ and $\tr E|_{\widetilde{X}} = \delta_{\mathbb{K}}(n)$. We define the desired operator $E$ as follows (recall that $C=\frac{\delta_{\mathbb{K}}(n)}{n}$):
\begin{enumerate}
\item $E(e_k)_j = \frac{\sgn \langle w_j, w_k \rangle}{d}$ for $1 \leq j, k \leq d$.
\item $E(e_k) = 0$ for $d+1 \leq k \leq N$.
\item $E(e_k)_{d+j} = C \overline{a_{jk}}$ for $1 \leq k \leq d$ and $1 \leq j \leq m$.
\end{enumerate}
Since $|Ca_{jk}| \leq \frac{1}{d}$ by the assumption, we clearly have
$$\sum_{i=1}^{N} \|E(e_i)\|_{\infty}=\sum_{i=1}^{d} \|E(e_i)\|_{\infty} = d \cdot \frac{1}{d}=1.$$
With a very similar calculation as before we shall verify that $E|_{\widetilde{X}}=  C\id_{\widetilde{X}}$. This will immediately imply that $E(\widetilde{X}) \subseteq \widetilde{X}$ and $\tr E|_{\widetilde{X}} = \delta_{\mathbb{K}}(n)$. Thus, let us fix $i \in \{1, \ldots, n\}$. If $1 \leq j \leq d$, then by Lemma \ref{lemwj}
\begin{align*}
E(x_i)_j =  \sum_{k=1}^{N} x_{ik} E(e_k)_j &= \sum_{k=1}^{d} x_{ik} E(e_k)_j = \frac{1}{d} \sum_{k=1}^{d} \overline{w_{ki} \sgn \langle w_k, w_j \rangle} \\
&=\frac{C}{d} \overline{w_{ji}} = C x_{ij}.
\end{align*}
On the other hand, if $1 \leq j \leq m$, then
\begin{align*}
E(x_i)_{d+j} =  \sum_{k=1}^{N} x_{ik} E(e_k)_{d+j} &= \sum_{k=1}^{d} x_{ik} E(e_k)_{d+j} = C \sum_{k=1}^{d} x_{ik} \overline{a_{jk}} \\
&=C \sum_{k=1}^{d} \overline{w_{ki} a_{jk}} = Cx_{i(d+j)}.
\end{align*}
It follows that the equality $E(x_i)=Cx_i$ holds for any $i \in \{1, \ldots, n\}$ and, consequently, $E|_{\widetilde{X}}= C \id_{\widetilde{X}}$. This concludes the proof. 

\end{proof}

\section{Characterization of spaces with the maximal projection constant in the general case}
\label{secchargeneral}

In this section we derive Theorem \ref{twrchar} in the general case by an approximation of a general $n$-dimensional normed space with a sequence of polyhedral spaces. Let us recall that for two $n$-dimensional normed spaces $X, Y$ over $\mathbb{K}$ their \emph{Banach--Mazur distance} is defined as
$$\dbm(X, Y) = \inf \|T\|\cdot \|T^{-1}\|,$$
where the infimum (in fact, a minimum) is taken over all invertible linear operators $T: X \to Y$ and the standard operator norm is considered. The set of isometric classes of $n$-dimensional normed spaces over $\mathbb{K}$ equipped with $\dbm$ is called the \emph{$n$-dimensional Banach--Mazur compactum}. It is a standard fact that 
$$ \frac{\lambda(X)}{\lambda(Y)}\le \dbm(X,Y)$$
for all $n$-dimensional normed spaces $X$ and $Y$ (see Corollary~6 in Section~III.B of \cite{wojtaszczyk}). Therefore, the mapping $X\mapsto\lambda(X)$ is continuous on the Banach--Mazur compactum. Consequently, for every sequence $(X_i)$ of polyhedral spaces satisfying $\dbm(X_i,X)\to 1$, we also have $\lambda(X_i)\to\lambda(X)$.

Using this observation, it is straightforward to deduce from the polyhedral case that any space $X$ satisfying the inclusions from Theorem~\ref{twrchar} has the maximal projection constant, since one can approximate $X$ by a sequence of polyhedral spaces attaining the maximal projection constant. However, proving the necessity of these inclusions is considerably more involved. In this situation, the approximating polyhedral spaces are only known to have almost maximal projection constants. Consequently, the arguments from Section~\ref{seccharpolyhedral} must be refined to handle an approximate equality. This requires precise asymptotic control of several quantities appearing in the proof of the polyhedral case, as now $N$ is no longer fixed but tends to infinity.

Our approach is to establish qualitative stability versions of the successive steps in the proof of the polyhedral case. We first develop a stability version of Theorem~\ref{twrbound}. To this end, we begin with a simple lemma.

\begin{lem} 
\label{operatornorm}
Let $(N_m)_{m=1}^{\infty}$ be a sequence of positive integers and let $ a^{(m)}, b^{(m)} \in \mathbb{R}^{N_m}$ be sequences of real vectors with non-negative coordinates such that 
$$ \|a^{(m)}-b^{(m)}\| \to 0 \qquad \text{ and } \qquad \|a^{(m)}\|,\ \|b^{(m)}\| \leq C, $$
for some constant $C>0$ independent of $m$. Suppose that $ G_m=(g_{ij}^{(m)})_{i,j=1}^{N_m} $ is a real matrix satisfying 
$$ 0\leq g_{ij}^{(m)} \leq a_i^{(m)}a_j^{(m)} $$ for every $i,j \in \{1, \ldots, N_m\}$. Then 
$$ \sum_{i,j=1}^{N_m} \Bigl( a_i^{(m)}a_j^{(m)} - b_i^{(m)}b_j^{(m)} \Bigr) g_{ij}^{(m)} \to 0. $$ 
\end{lem} 
\begin{proof} For simplicity, we omit the superscript $(m)$ and the summation will be always over $\{1, \ldots, N_m\}$. Since 
$$ a_i a_j-b_i b_j = a_i(a_j-b_j)+b_j(a_i-b_i), $$ 
we obtain 
$$ |a_i a_j-b_i b_j| \leq a_i|a_j-b_j| + b_j|a_i-b_i|. $$ Hence 
$$\sum_{i,j} |a_i a_j-b_i b_j| g_{ij}\leq \sum_{i,j} \Bigl( a_i|a_j-b_j| + b_j|a_i-b_i| \Bigr) g_{ij}.$$
Using the bound $ g_{ij}\leq a_i a_j$, we get
$$\sum_{i,j} |a_i a_j-b_i b_j| g_{ij} \leq \sum_{i,j} a_i^2 a_j |a_j-b_j| + \sum_{i,j} a_i a_j b_j |a_i-b_i|.$$
Factoring the sums gives 
\begin{align*} \sum_{i,j} a_i^2 a_j |a_j-b_j| &= \Bigl(\sum_i a_i^2\Bigr) \Bigl(\sum_j a_j|a_j-b_j|\Bigr), \\ \sum_{i,j} a_i a_j b_j |a_i-b_i| &= \Bigl(\sum_j a_j b_j\Bigr) \Bigl(\sum_i a_i|a_i-b_i|\Bigr). \end{align*}
By the Cauchy--Schwarz inequality, 
$$ \sum_i a_i|a_i-b_i| \leq \|a\|\,\|a-b\|, $$ and similarly 
$$ \sum_j a_j|a_j-b_j| \leq \|a\|\,\|a-b\|, $$ while 
$$ \sum_j a_j b_j \leq \|a\|\,\|b\|. $$ 
Consequently
$$\sum_{i,j} |a_i a_j-b_i b_j| g_{ij} \leq \|a\|^3\,\|a-b\| + \|a\|^2\|b\|\,\|a-b\|.$$
From $ \|a^{(m)}\|,\ \|b^{(m)}\|\leq C$ it now follows that 
$$ \sum_{i,j} |a_i^{(m)}a_j^{(m)}-b_i^{(m)}b_j^{(m)}| g_{ij}^{(m)} \leq 2C^3\, \|a^{(m)}-b^{(m)}\|. $$ 
Since $ \|a^{(m)}-b^{(m)}\| \to 0$, the right-hand side converges to $0$, proving the lemma. 
\end{proof}

Following the terminology of \cite{bukhcox}, an \emph{isotropic measure} on $\mathbb{K}^n$ is a probabilistic measure $\mu$ such that $\int_{\mathbb{K}^n} xx^* \,d\mu(x) = \frac{1}{n} \id_n$. This generalizes the concept of tight frames, since for a tight frame $u_1, \ldots, u_N \in \mathbb{K}^n$ with constant $1$, the measure $\mu$ defined by
$$\mu = \sum_{i=1}^{N} \frac{\|u_i\|^2}{n} \delta_{u_i/\|u_i\|},$$
(where $\delta_x$ denotes the Dirac measure at $x$) is isotropic. In Lemma 5 of \cite{bukhcox} it was proved that $\iint |\langle x,y\rangle| \,d\mu(x)\,d\mu(y) \leq \frac{\delta_{\mathbb{K}}(n)}{n}$ for every isotropic measure $\mu$ on $\mathbb{K}^n$. An isotropic measure will be called an \emph{extremal isotropic measure} if equality holds in this estimate. Lemma~5 of \cite{bukhcox} characterizes the extremal isotropic measures by conditions similar to those appearing in Theorem~\ref{twrbound}: for equality to hold, a maximal ETF must exist in $\mathbb K^n$, and the measure must be equidistributed among all ETF directions, each receiving measure $\frac1d$. It should be noted that part~(b) of Lemma~5 in \cite{bukhcox} requires a correction in the complex case. As pointed out in the published corrigendum \cite{bukhcoxcorr}, the condition should state that the measure of the entire complex circle $\{\alpha x:\ |\alpha|=1\}$ (where $x$ is an ETF direction) is equal to $\frac1d$, rather than merely the total measure of the vectors $x$ and $-x$. The proof of Lemma~5 from \cite{bukhcox} remains valid after this modification. Let us also point out that the result of Bukh and Cox concerns only isotropic measures and does not involve weights $t_i$ appearing in \eqref{formula}. However, the proof of Theorem~\ref{twrbound} makes it possible to extract asymptotic information about the weights $t_i$ in the case of an asymptotic equality. To apply the characterization of extremal isotropic measures due to Bukh and Cox, we shall use weak* compactness of measures together with Prokhorov's theorem to extract a weak* convergent subsequence from the sequence of measures induced by the tight frame vectors.

\begin{twr}
\label{twrboundstab}
Let $n \geq 2$ be an integer such that there exists a maximal ETF in $\mathbb{K}^n$ (of cardinality $d = d_{\mathbb{K}}(n)$). For every $m \geq 1$ suppose that $N_m \geq n$ is an integer and vectors $u_1^{(m)}, \ldots, u^{(m)}_{N_m} \in \mathbb{K}^n$ form a tight frame in $\mathbb{K}^n$ with constant $1$. Let $t^{(m)} \in \mathbb{R}^{N_m}$ be a vector with non-negative coordinates satisfying $\|t^{(m)}\|=1$. Assume that
$$\delta_m=\sum_{i, j=1}^{N_m} t^{(m)}_it^{(m)}_j |\langle u^{(m)}_i, u^{(m)}_j \rangle| \to \delta_{\mathbb{K}}(n),$$
as $m \to \infty$. Then there exists a maximal ETF $w_1, \ldots, w_{d}$ in $\mathbb{K}^n$ and a subsequence (which we do not relabel) of the considered vectors, along with partitions $\{1, \ldots, N_m\} = A^{(m)}_0 \sqcup A^{(m)}_1 \sqcup \ldots \sqcup A^{(m)}_d$ satisfying the following conditions:
\begin{enumerate}
\item $\sum_{i \in A^{(m)}_0} \left (t_i^{(m)} \right)^2 \to 0$ and $\sum_{i \in A^{(m)}_0} \left \|u_i^{(m)} \right \|^2 \to 0$. Moreover, the set $A^{(m)}_0$ contains all indices $i$ such that $u_i^{(m)}=0$.
\item For every $j=1, \ldots, d$ we have
$$\sum_{i \in A^{(m)}_j} (t^{(m)}_i)^2 \inf_{|\alpha|=1} \left \| \frac{u_i^{(m)}}{\|u_i^{(m)}\|} - \alpha w_j \right \|^2 \to 0.$$
\item For every $j=1, \ldots, d$ we have $\sum_{i \in A^{(m)}_j} \left (t^{(m)}_i \right )^2 \to \frac{1}{d}$.
\item $\sum_{i=1}^{N_m} \left ( \left \| u_i^{(m)} \right  \| - \sqrt{n}t_i^{(m)} \right )^2 \to 0.$
\end{enumerate}
\end{twr}

\begin{proof}
Let $S^{(m)} \subseteq \{1, \ldots, N_m\}$ be the set of indices $i$ such that $u^{(m)}_i=0$ and $I^{(m)}=\{1, \ldots, N_m\} \setminus S^{(m)}$. We shall first justify that the condition $(4)$ holds, even without passing to a subsequence. Clearly, the vectors $u^{(m)}_i$ for $i \in I^{(m)}$ still form a tight frame with constant $1$, with the associated vector $(t')^{(m)} \in \mathbb{R}^{|I^{(m)}|}$ obtained by removing the coordinates corresponding to indices from $S^{(m)}$. Repeating the exact same estimates as in the proof of Theorem \ref{twrbound} we arrive at the estimate analogous to \eqref{cst}:
$$\left( \sum_{i \in I^{(m)}} (t_i')^{(m)} \|u^{(m)}_i\| \right)^2 \leq \|(t')^{(m)}\|^2 \cdot \left ( \sum_{i \in I^{(m)}} \|u^{(m)}_i\|^2 \right ) = \|(t')^{(m)}\|^2 \cdot n \leq n,$$

from which later it follows that 
$$\sum_{i, j \in I^{(m)}} (t')^{(m)}_i (t')^{(m)}_j |\langle u^{(m)}_i, u^{(m)}_j \rangle| \leq \frac{1}{2\varphi} + \frac{1}{2 \varphi} \left ( \varphi^2 - \frac{(1-\varphi)^2}{d} \right) \left( \sum_{i \in I^{(m)}} (t')^{(m)}_i \|u^{(m)}_i\| \right)^2$$
$$\leq \frac{1}{2 \varphi} + \frac{n}{2 \varphi} \left ( \varphi^2 - \frac{(1-\varphi)^2}{d} \right) = \delta_{\mathbb{K}}(n).$$
Since 
$$ \sum_{i,j=1}^{N_m} t_i^{(m)}t_j^{(m)} |\langle u_i^{(m)},u_j^{(m)}\rangle| \to\delta_{\mathbb K}(n),$$
all inequalities appearing above are asymptotically sharp. Consequently, we obtain $\sum_{i \in I^{(m)}} \left (t^{(m)}_i \right )^2 \to 1$ (or equivalently $\sum_{i \in S^{(m)}} \left ( t^{(m)}_i \right )^2 \to 0$) and 
$$ \sum_{i=1}^{N_m} t^{(m)}_i \|u^{(m)}_i\|  \to \sqrt{n}.$$
Hence
$$\sum_{i=1}^{N_m} \left ( \left \| u_i^{(m)} \right  \| - \sqrt{n}t_i^{(m)} \right )^2 = \sum_{i=1}^{N_m} \left \| u_i^{(m)} \right  \|^2 + n \sum_{i=1}^{N_m} \left ( t_i^{(m)} \right  )^2  - 2 \sqrt{n}\sum_{i=1}^{N_m} t^{(m)}_i \|u^{(m)}_i\|$$
$$=2n - 2\sqrt{n} \sum_{i=1}^{N_m} t^{(m)}_i \|u^{(m)}_i\| \to 0,$$
as claimed.

For every $m\geq 1$ define a measure on the unit sphere $\mathbb{S}^{n-1}$ of $\mathbb K^n$ by 
$$ \nu_m = \sum_{i\in I^{(m)}} \frac{\|u_i^{(m)}\|^2}{n}\, \delta_{u_i^{(m)}/\|u_i^{(m)}\|}.$$ 
Since the vectors $u_i^{(m)}$ form a tight frame with constant $1$, we have $ \sum_{i=1}^{N_m}\|u_i^{(m)}\|^2=n$ and thus 
$$ \nu_m(\mathbb{S}^{n-1}) = \sum_{i\in I^{(m)}}\frac{\|u_i^{(m)}\|^2}{n} = 1, $$
so each $\nu_m$ is a probability measure. Moreover, by the tight frame identity 
$$ \sum_{i\in I^{(m)}}u_i^{(m)}(u_i^{(m)})^*=\id_n, $$ 
we obtain 
$$ n\int_{\mathbb{S}^{n-1}}xx^*\,d\nu_m(x) = n\sum_{i\in I^{(m)}} \frac{\|u_i^{(m)}\|^2}{n} \frac{u_i^{(m)}}{\|u_i^{(m)}\|} \left( \frac{u_i^{(m)}}{\|u_i^{(m)}\|} \right)^* = \id_n. $$ Hence every measure $\nu_m$ is isotropic. Since the unit sphere $\mathbb{S}^{n-1}$ is compact, the sequence $(\nu_m)$ is tight. By Prokhorov's theorem, passing to a subsequence if necessary, we may assume that 
$$ \nu_m \stackrel{*}{\to} \nu $$ 
for some Borel probability measure $\nu$ on the sphere. Since the map 
$$ \mathbb{S}^{n-1}\ni x\mapsto xx^*\in \mathbb{K}^{n^2} $$ 
is continuous, weak* convergence yields 
$$ \int_{\mathbb{S}^{n-1}}xx^*\,d\nu_m(x) \to \int_{\mathbb{S}^{n-1}}xx^*\,d\nu(x). $$ 
Passing to the limit in 
$$ n\int_{\mathbb{S}^{n-1}}xx^*\,d\nu_m(x)=\id_n $$ 
gives $$ n\int_{\mathbb{S}^{n-1}}xx^*\,d\nu(x)=\id_n. $$ 
Thus $\nu$ is isotropic. Let 
$$ F(\mu) = \iint |\langle x,y\rangle| \,d\mu(x)\,d\mu(y). $$ 
Since the mapping $(x,y)\mapsto |\langle x,y\rangle|$ is continuous and bounded on $\mathbb{S}^{n-1} \times \mathbb{S}^{n-1}$, and weak* convergence $ \nu_m \stackrel{*}{\to} \nu $ implies weak* convergence of the product measures, we have $F(\nu_m)\to F(\nu)$. Furthermore, 
\begin{align*} F(\nu_m) &= \sum_{i,j\in I^{(m)}} \frac{\|u_i^{(m)}\|^2}{n} \frac{\|u_j^{(m)}\|^2}{n} \left| \left\langle \frac{u_i^{(m)}}{\|u_i^{(m)}\|}, \frac{u_j^{(m)}}{\|u_j^{(m)}\|} \right\rangle \right|\\ &= \frac1{n^2} \sum_{i,j\in I^{(m)}} \|u_i^{(m)}\| \|u_j^{(m)}\| |\langle u_i^{(m)},u_j^{(m)}\rangle|. \end{align*} 
Now, let us take 
$$ a^{(m)} = \bigl(\|u_1^{(m)}\|,\ldots,\|u_{N_m}^{(m)}\|\bigr) \in \ell_2^{N_m} \quad \text{ and } \quad b^{(m)} = \bigl(\sqrt n\,t_1^{(m)},\ldots,\sqrt n\,t_{N_m}^{(m)}\bigr) \in \ell_2^{N_m}. $$ 
By the already derived condition (4) we have $ \|a^{(m)}-b^{(m)}\| \to 0. $
Moreover, $ \|a^{(m)}\|^2 = \|b^{(m)}\|^2 = n$. Lemma \ref{operatornorm} applied to the vectors $a^{(m)}$, $b^{(m)}$ and the matrix $$ G_m= \left (|\langle u_i^{(m)},u_j^{(m)}\rangle| \right )_{i,j=1}^{N_m}$$
yields that
$$\left | F(\nu_m) - \frac1n \sum_{i,j=1}^{N_m} t_i^{(m)}t_j^{(m)} |\langle u_i^{(m)},u_j^{(m)}\rangle| \right | = \left | F(\nu_m) - \frac{\delta_m}{n} \right | \to 0.$$
Hence, passing to the limit leads to
$$ F(\nu) = \frac{\delta_{\mathbb K}(n)}{n}.$$ 
It follows that $\nu$ is an extremal isotropic measure. Let $ \widetilde\nu = \pi_\#(\nu) $ be the push-forward of $\nu$ to the projective space under the quotient map $ \pi:\mathbb S^{n-1}\to\mathbb P(\mathbb K^n)$. From the characterization of extremal isotropic measures given in Lemma 5 of \cite{bukhcox}, it follows that there exists a maximal ETF $ w_1,\ldots,w_d$ in  $\mathbb{K}^n $ such that
$$ \widetilde\nu = \frac1d \sum_{j=1}^{d}\delta_{[w_j]}. $$ 
Since the projective space $\mathbb P(\mathbb K^n)$ is compact, we may choose open neighborhoods $ U_1,\ldots,U_d $ of the projective points $ [w_1],\ldots,[w_d]$ with pairwise disjoint closures (in particular $\widetilde\nu(\partial U_j)=0$ for every $j$). For every $m\geq 1$ and $j \in \{1, \ldots, d\}$ we define 
$$ A_j^{(m)} = \left\{ i\in I^{(m)}: \ \pi\!\left( \frac{u_i^{(m)}}{\|u_i^{(m)}\|} \right)\in U_j \right\}$$ 
and 
$$ A_0^{(m)} = \{1,\ldots,N_m\} \setminus \bigcup_{j=1}^{d}A_j^{(m)}.$$
Since $\widetilde{\nu}=\frac1d\sum_{j=1}^{d}\delta_{[w_j]}$, weak* convergence of the measures $\nu_m$ to $\nu$ yields the corresponding convergence for push-forward measures $\widetilde{\nu_m}$ and hence
$$ \sum_{i\in A_j^{(m)}} \frac{\|u_i^{(m)}\|^2}{n} = \widetilde{\nu_m}(U_j) \to \widetilde{\nu}(U_j) = \frac1d, $$
by Portmanteau theorem, as $\widetilde\nu(\partial U_j)=0$ for every $j$. We now claim that
$$ \sum_{i\in A_j^{(m)}} \bigl(t_i^{(m)}\bigr)^2 \to \frac1d, $$ 
which will establish condition $(3)$. Indeed, it follows from the Cauchy--Schwarz inequality, condition $(4)$ and the equalities $\sum_{i=1}^{N_m} \|u_i^{(m)}\|^2 = \sum_{i=1}^{N_m} n(t_i^{(m)})^2=n$ that
\begin{align*}
\sum_{i=1}^{N_m} \left| \|u_i^{(m)}\|^2 - n(t_i^{(m)})^2 \right| =& \sum_{i=1}^{N_m} \left| \|u_i^{(m)}\|-\sqrt n\,t_i^{(m)} \right| \left| \|u_i^{(m)}\|+\sqrt n\,t_i^{(m)} \right| \\
\leq& \sqrt{\sum_{i=1}^{N_m} \Bigl( \|u_i^{(m)}\|-\sqrt n\,t_i^{(m)} \Bigr)^2 \sum_{i=1}^{N_m} \left ( \|u_i^{(m)}\|+\sqrt n\,t_i^{(m)} \right )^2 } \\
\leq& 2\sqrt{n} \sqrt{\sum_{i=1}^{N_m} \Bigl( \|u_i^{(m)}\|-\sqrt n\,t_i^{(m)} \Bigr)^2} \to 0.
\end{align*}
Hence
$$\left | \sum_{i\in A_j^{(m)}} \frac{\|u_i^{(m)}\|^2}{n} - \sum_{i\in A_j^{(m)}} (t_i^{(m)})^2 \right | \leq \frac{1}{n} \sum_{i\in A_j^{(m)}} \left| \|u_i^{(m)}\|^2 - n(t_i^{(m)})^2 \right| \to 0 ,$$
which yields the desired convergence.  Moreover, 
$$ \widetilde{\nu}\!\left( \mathbb P(\mathbb K^n) \setminus \bigcup_{j=1}^{d}U_j \right) =0, $$
and therefore 
$$ \sum_{i\in A_0^{(m)}} \frac{\|u_i^{(m)}\|^2}{n} = \widetilde{\nu_m}\!\left( \mathbb P(\mathbb K^n) \setminus \bigcup_{j=1}^{d}U_j \right) \to 0. $$ Hence $ \sum_{i\in A_0^{(m)}}\|u_i^{(m)}\|^2\to0, $ and similarly as before we obtain that $ \sum_{i\in A_0^{(m)}} \bigl(t_i^{(m)}\bigr)^2 \to 0$. This proves condition (1).

For a fixed $j = 1,\ldots,d$ let us now define
$$ d_i^{(m)} = \inf_{|\alpha|=1} \left\| \frac{u_i^{(m)}}{\|u_i^{(m)}\|} - \alpha w_j \right\|, \qquad \text{ for } i\in A_j^{(m)}. $$  
To establish condition $(2)$ we shall first prove that $a_m \to 0$, where
$$ a_m = \sum_{i\in A_j^{(m)}} \frac{\|u_i^{(m)}\|^2}{n} \left(d_i^{(m)}\right)^2. $$ Fix $\varepsilon>0$ and define 
$$ B_\varepsilon^{(m)} = \left\{ i\in A_j^{(m)}: d_i^{(m)}\geq \varepsilon \right\}. $$ Then 
$$a_m = \sum_{i\in A_j^{(m)}\setminus B_\varepsilon^{(m)}} \frac{\|u_i^{(m)}\|^2}{n} \left(d_i^{(m)}\right)^2 + \sum_{i\in B_\varepsilon^{(m)}} \frac{\|u_i^{(m)}\|^2}{n} \left(d_i^{(m)}\right)^2.$$ 
Since $d_i^{(m)}<\varepsilon$ for $i\in A_j^{(m)}\setminus B_\varepsilon^{(m)}$, we obtain 
\begin{align*} \sum_{i\in A_j^{(m)}\setminus B_\varepsilon^{(m)}} \frac{\|u_i^{(m)}\|^2}{n} \left(d_i^{(m)}\right)^2 &\leq \varepsilon^2 \sum_{i\in A_j^{(m)}} \frac{\|u_i^{(m)}\|^2}{n}. \end{align*}
From $$ \sum_{i\in A_j^{(m)}} \frac{\|u_i^{(m)}\|^2}{n} \to \frac1d, $$ it follows now that 
$$ \limsup_{m\to\infty} \sum_{i\in A_j^{(m)}\setminus B_\varepsilon^{(m)}} \frac{\|u_i^{(m)}\|^2}{n} \left(d_i^{(m)}\right)^2 \leq \frac{\varepsilon^2}{d}.$$
Now let 
$$ F_\varepsilon = \Bigl\{ [x]\in \overline{U}_j: \inf_{|\alpha|=1} \left\| x - \alpha w_j \right\|\geq \varepsilon \Bigr\},$$
where $x$ denotes any unit representative of the projective point $[x]$. Since $ \widetilde{\nu}=\frac1d\sum_{r=1}^{d}\delta_{[w_r]}$, we have $ \widetilde{\nu}(F_\varepsilon)=0$. Thus, weak* convergence and Portmanteau theorem give $ \widetilde{\nu_m}(F_\varepsilon) \to 0$, as $F_{\varepsilon}$ is a closed set of $\widetilde{\nu}$-measure $0$. By definition $ \widetilde{\nu_m}(F_\varepsilon) = \sum_{i\in B_\varepsilon^{(m)}} \frac{\|u_i^{(m)}\|^2}{n}$ and hence
$$ \sum_{i\in B_\varepsilon^{(m)}} \frac{\|u_i^{(m)}\|^2}{n} \to 0. $$
Since $d_i^{(m)}\leq 2$, we obtain 
\begin{align*} \sum_{i\in B_\varepsilon^{(m)}} \frac{\|u_i^{(m)}\|^2}{n} \left(d_i^{(m)}\right)^2 &\leq 4 \sum_{i\in B_\varepsilon^{(m)}} \frac{\|u_i^{(m)}\|^2}{n} \to 0. \end{align*} 
Combining the above estimates yields $ \limsup_{m\to\infty} a_m \leq \frac{\varepsilon^2}{d}.$  Since $\varepsilon>0$ was arbitrary, we conclude that $ \limsup_{m\to\infty} a_m=0$, so that $a_m \to 0$. Now by $d_i^{(m)}\leq 2$ we get
\begin{align*}
\left |\sum_{i\in A_j^{(m)}} \left ( \frac{\|u_i^{(m)}\|^2}{n} - (t_i^{(m)})^2 \right ) \left(d_i^{(m)}\right)^2 \right | &\leq 4 \sum_{i\in A_j^{(m)}} \left| \|u_i^{(m)}\|^2 - n(t_i^{(m)})^2 \right| \\
&\leq 4 \sum_{i=1}^{N_m} \left| \|u_i^{(m)}\|^2 - n(t_i^{(m)})^2 \right| \to 0,
\end{align*}
as before. This establishes condition $(2)$ and concludes the proof.

\end{proof}

The following technical lemma is a stability version of Lemma~\ref{lemtiui}.

\begin{lem} 
\label{lemtiuistab} 
Let $n\geq 2$ be an integer such that there exists a maximal ETF in $\mathbb K^n$. 
Let $(N_m)_{m=1}^\infty$ be a sequence of integers with $N_m\geq n$. Suppose that vectors $u_1^{(m)},\ldots,u_{N_m}^{(m)}\in\mathbb K^n$ form a tight frame with constant $1$ and let $t^{(m)}\in\mathbb R^{N_m}$ be a vector with non-negative coordinates satisfying $\|t^{(m)}\|=1$. Assume that 
$$ \sum_{i,j=1}^{N_m} t_i^{(m)}t_j^{(m)} |\langle u_i^{(m)},u_j^{(m)}\rangle| \to \delta_{\mathbb K}(n). $$ 
Then there exists a subsequence (which we do not relabel) such that
$$\sum_{j=1}^{N_m} \left\| u_j^{(m)} - \frac{nt_j^{(m)}}{\delta_{\mathbb K}(n)} \sum_{i=1}^{N_m} t_i^{(m)} \sgn \langle u_i^{(m)}, u_j^{(m)}\rangle u_i^{(m)} \right\|^2 \to 0. $$ 
\end{lem}

\begin{proof} Throughout the proof, for clarity we drop the superscript $(m)$ from all sequences, understanding that limits are taken as $m \to \infty$. We denote $d = d_{\mathbb{K}}(n)$, $\delta = \delta_{\mathbb{K}}(n)$ and $\varphi=\varphi_{\mathbb{K}}(n)$. For $i, j \in \{1, \ldots, N\}$ let also $H_{ij} = \sgn\langle u_i, u_j \rangle$ and  $v_j = \sum_{i=1}^N t_i H_{ij} u_i$. Our goal is to prove that $\sum_{j=1}^N \left \| u_j - \frac{n}{\delta} t_j v_j \right \|^2 \to 0$.

Let us pass to the subsequence for which the partition
$$ \{1,\ldots,N\} = A_0 \sqcup A_1 \sqcup \ldots \sqcup A_d $$ 
and the ETF vectors $ w_1,\ldots,w_d$ are given by Theorem~\ref{twrboundstab}. For $i \in \{1, \ldots N\} \setminus A_0$ let $s(i) \in \{1, \ldots, d\}$ be unique number such that $i \in A_{s(i)}$. Moreover, for every such $i$ we choose a unimodular scalar $\alpha_i$ satisfying 
$$ d_i = \left\| \frac{u_i}{\|u_i\|} - \alpha_iw_{s(i)} \right\| = \inf_{|\alpha|=1} \left\| \frac{u_i}{\|u_i\|} - \alpha w_{s(i)} \right\|. $$ 
We further define $\widetilde{u}_i = \sqrt{n} t_i \alpha_i w_{s(i)}$ for $i \notin A_0$. We first deal with the exceptional set $A_0$. Condition $(1)$ of Theorem~\ref{twrboundstab} implies
$$ \sum_{i\in A_0}\|u_i\|^2\to 0 \qquad \text{ and } \qquad \sum_{i\in A_0}t_i^2\to 0. $$ 
Moreover, for any $j=1, \ldots, N$ we have
$$\|v_j\|=\left\| \sum_{k=1}^{N} t_k \sgn \langle u_k,u_j\rangle u_k \right\| \leq \sum_{k=1}^{N} t_k\|u_k\| \leq \sqrt{\Bigl( \sum_{k=1}^{N} (t_k)^2 \Bigr) \Bigl( \sum_{k=1}^{N} \|u_k\|^2 \Bigr)} = \sqrt n.$$
Therefore 
$$ \sum_{j\in A_0}\left \|u_j - \frac{n}{\delta} t_j v_j \right \|^2 \leq 2 \left (\sum_{j \in A_0 }\|u_j\|^2 +  \frac{n^3}{\delta^2_{\mathbb K}(n)} \sum_{j \in A_0} t_j^2 \right ) \to 0. $$
Thus, from now on we shall work only with indices from $\{1, \ldots N\} \setminus A_0$. We first establish that 
\begin{equation}
\label{conv2}
\sum_{i \notin A_0} \|u_i - \widetilde{u}_i\|^2 \to 0.
\end{equation}
Indeed, for $i \notin A_0$ we have
$$ u_i - \widetilde{u}_i = \Bigl( \|u_i\| - \sqrt n\,t_i \Bigr) \frac{u_i}{\|u_i\|} + \sqrt n\,t_i \left( \frac{u_i}{\|u_i\|} - \alpha_iw_{s(i)} \right).$$
Hence 
$$\left\| u_i - \widetilde{u}_i \right\|^2 \leq 2 \Bigl( \|u_i\| - \sqrt n\,t_i \Bigr)^2 + 2n (t_id_i)^2.$$
Summing over all $i\notin A_0$ and using conditions $(2)$ and $(4)$ of Theorem~\ref{twrboundstab} yields \eqref{conv2}.

For $j \not \in A_0$ let us now define 
$$\widetilde{v}_j = \sum_{i \notin A_0} t_i \widetilde{H}_{ij} \widetilde{u}_i,$$
where 
$$\widetilde{H}_{ij} = \sgn\langle \widetilde{u}_i, \widetilde{u}_j \rangle = \overline{\alpha_i} \alpha_j \sgn \langle w_{s(i)}, w_{s(j)} \rangle.$$
Since $\overline{\alpha}_i \alpha_i = 1$ we have
$$\widetilde{H}_{ij} \widetilde{u}_i = \sqrt{n} t_i \alpha_j  \sgn \langle w_{s(i)}, w_{s(j)} \rangle w_{s(i)}.$$
Substituting this into $\widetilde{v}_j$ and grouping by the partition sets we get:
$$\widetilde{v}_j = \sqrt{n} \alpha_j \sum_{s=1}^d \sgn \langle w_{s}, w_{s(j)} \rangle w_s \left( \sum_{i \in A_s} t_i^2 \right).$$
By condition $(3)$ of Theorem~\ref{twrboundstab}, we have $\sum_{i \in A_s} t_i^2 = \frac{1}{d} + \gamma_s$ where $\gamma_s \to 0$. From Lemma \ref{lemtiui} it follows that
$$\widetilde{v}_j = \frac{\delta}{\sqrt{n}} \alpha_j w_{s(j)} + R_j = \frac{ \delta \widetilde{u}_j}{n t_j} + R_j,$$
where $R_j = \sqrt{n} \alpha_j \sum_{s=1}^d \gamma_s \sgn \langle w_{s}, w_{s(j)} \rangle w_s$. Since $|\gamma_s| \to 0$, we have $\max_{j \notin A_0} \|R_j\| \to 0$.
Hence
\begin{equation}
\label{conv3}
\sum_{j \notin A_0} \left \| \widetilde{u}_j - \frac{n}{\delta}t_j \widetilde{v}_j \right \|^2 = \frac{n^2}{\delta^2} \sum_{j \notin A_0} t_j^2 \|R_j\|^2 \le \frac{n^2}{\delta^2}\left( \max_{j \notin A_0} \|R_j\|^2 \right) \to 0.
\end{equation}
For $j \notin A_0$ we now take $x_j = t_j v_j$, $\widetilde{x}_j = t_j \widetilde{v}_j$ and consider the difference
$$x_j - \widetilde{x}_j = \underbrace{\sum_{i \notin A_0} t_i t_j H_{ij} (u_i - \widetilde{u}_i)}_{:= A_j} + \underbrace{\sum_{i \notin A_0} t_i t_j (H_{ij} - \widetilde{H}_{ij}) \widetilde{u}_i}_{:= B_j} + \underbrace{\sum_{i \in A_0} t_i t_j H_{ij}u_i}_{:= C_j}$$
We note that to prove the lemma it suffices to establish that $\sum_{j \notin A_0} \|x_j-\widetilde{x}_j\|^2 \to 0$, because
$$ \left \|u_j - \frac{n}{\delta} t_j v_j \right \| = \left \|u_j - \frac{n}{\delta} x_j \right \| \leq \|u_j - \widetilde{u}_j \| + \left \| \widetilde{u}_j - \frac{n}{\delta} \widetilde{x}_j \right \| + \frac{n}{\delta}\left \|x_j-\widetilde{x}_j \right \|, $$
and the sums of squares the first two terms on the right-hand side go to $0$ by \eqref{conv2} and \eqref{conv3}. Thus, in the remaining part of the proof we establish that $\sum_{j \notin A_0} \|A_j\|^2 \to 0$, $\sum_{j \notin A_0} \|B_j\|^2 \to 0$ and $\sum_{j \notin A_0} \|C_j\|^2 \to 0$.

The first and third convergences are immediate, since by the Cauchy--Schwarz inequality we can estimate 
$$\|A_j\|^2 \le t_j^2 \left(\sum_{i \notin A_0} t_i^2 \right ) \left(  \sum_{i \notin A_0} \|u_i - \widetilde{u}_i\|^2 \right ) \leq t_j^2\sum_{i \notin A_0} \|u_i - \widetilde{u}_i\|^2.$$
Summing over $j$ gives 
$$\sum_{j \not \in A_0} \|A_j\|^2 \leq \left ( \sum_{j \notin A_0} t_j^2  \right ) \sum_{i \notin A_0} \|u_i - \widetilde{u}_i\|^2 \leq \sum_{i \notin A_0} \|u_i - \widetilde{u}_i\|^2 \to 0$$
by \eqref{conv2}. Similarly, for $C_j$, by $|H_{ij}| \leq 1$ we can estimate
$$\|C_j\|^2 \leq t_j^2 \left ( \sum_{i \in A_0} t_i^2 \right )\left ( \sum_{i \in A_0} \|u_i\|^2 \right ) \leq \sum_{i \in A_0} \|u_i\|^2,$$
so that 
$$\sum_{j \notin A_0} \|C_j\|^2 \leq \left (\sum_{j \notin A_0} t_j^2 \right )\left ( \sum_{i \in A_0} \|u_i\|^2 \right ) \leq \sum_{i \in A_0} \|u_i\|^2 \to 0,$$
by condition (1) of Theorem \ref{twrboundstab}. For $B_j$, the Cauchy--Schwarz inequality again yields 
$$\|B_j\|^2 \le \left ( \sum_{i \notin A_0} \|\widetilde{u}_i\|^2 \right ) \left ( \sum_{i \notin A_0} t_i^2 t_j^2 | \Delta_{ij}|^2 \right ) \le n \sum_{i \notin A_0} t_i^2 t_j^2 |\Delta_{ij}|^2,$$
where $\Delta_{ij} = H_{ij} - \widetilde{H}_{ij}$. Summing over $j$ leads to
\begin{equation}
\label{eqbj}
\sum_{j \not\in A_0} \|B_j\|^2 \leq n \sum_{i,j \notin A_0} t_i^2 t_j^2 |\Delta_{ij}|^2.
\end{equation}
In order to prove that the quantity on the right-hand side goes to $0$ , we will use the identity
$$|z|-\re(z w) = \frac{1}{2} |z| | \sgn(z) - w|^2,$$
that holds for all $z, w \in \mathbb{K}$ with $|w|=1$. In particular, we have
$$\left | \langle \widetilde{u}_i, \widetilde{u}_j \rangle \right |- \re \left ( \langle \widetilde{u}_i, \widetilde{u}_j \rangle H_{ij} \right ) = \frac{1}{2} |\langle \widetilde{u}_i, \widetilde{u}_j \rangle| |\Delta_{ij}|^2.$$
Summing this weighted by $t_i t_j$ over all $i,j \notin A_0$ yields
$$\frac{1}{2} \sum_{i,j \notin A_0} t_i t_j |\langle \widetilde{u}_i, \widetilde{u}_j \rangle| |\Delta_{ij}|^2 = \widetilde{\delta}_m - \re \left ( \sum_{i,j \notin A_0} t_i t_j \langle \widetilde{u}_i, \widetilde{u}_j \rangle H_{ij} \right ).$$
where
$$\widetilde{\delta}_m = \sum_{i,j \notin A_0} t_i t_j |\langle \widetilde{u}_i, \widetilde{u}_j \rangle|.$$
By taking $E_{ij} = \langle u_i, u_j \rangle - \langle \widetilde{u}_i, \widetilde{u}_j \rangle $ for $i, j \notin A_0$ and noting that $\re(\langle u_i, u_j \rangle H_{ij}) = |\langle u_i, u_j \rangle |$ we obtain
\begin{equation}
\label{eqdeltas}
\frac{1}{2} \sum_{i,j \notin A_0} t_i t_j |\langle \widetilde{u}_i, \widetilde{u}_j \rangle| |\Delta_{ij}|^2 = \widetilde{\delta}_m - \delta'_m + \re \left ( \sum_{i, j \notin A_0} t_i t_j E_{ij} H_{ij} \right ),
\end{equation}
where
\begin{align*}\delta'_m = \sum_{i, j \notin A_0} t_it_j |\langle u_i, u_j \rangle| &= \sum_{i, j=1}^{N} t_it_j |\langle u_i, u_j \rangle| - \sum_{i \in A_0 \text{ or } j \in A_0} t_it_j |\langle u_i, u_j \rangle| \\
&=\delta_m - \sum_{i \in A_0 \text{ or } j \in A_0} t_it_j |\langle u_i, u_j \rangle|.
\end{align*}
We note that the subtracted sum disappears in the limit. Indeed
$$\sum_{i \in \{1, \ldots, N\}, j \in A_0} t_it_j|\langle u_i, u_j \rangle| \leq \sum_{i \in \{1, \ldots, N\}, j \in A_0} t_it_j\|u_i\| \|u_j\| = \left ( \sum_{i=1}^{N} t_i \|u_i\| \right ) \left ( \sum_{j \in A_0} t_j \|u_j\| \right ),$$
and by the Cauchy--Schwarz inequality the first factor remains bounded, while the second one goes to $0$ by condition (1) of Theorem \ref{twrboundstab}. The same argument applies with the roles of $i$ and $j$ switched, and hence $\delta'_m \to \delta_m$. Furthermore, since the vectors $\widetilde{u}_i$ are just rescaled ETF vectors, a similar computation like in the last part of the proof of Theorem \ref{twrbound} shows that
$$\widetilde{\delta}_m = 2n \varphi \sum_{1 \leq k < l \leq d} \left ( \sum_{i \in A_k} t_i^2 \right ) \left ( \sum_{j \in A_l} t_j^2 \right ) + n \sum_{k=1}^{d}\left ( \sum_{i \in A_k} t_i^2 \right )^2 \to \delta,$$
by condition $(3)$ of Theorem \ref{twrboundstab}. In particular, $\widetilde{\delta}_m - \delta'_m \to \delta - \delta = 0$. To show that the right-hand side of \eqref{eqdeltas} goes to $0$, it is therefore enough to establish that $\sum t_i t_j |E_{ij}| \to 0$. Indeed, $E_{ij} = \langle u_i, u_j \rangle - \langle \widetilde{u}_i, \widetilde{u}_j \rangle $, so
$$|E_{ij}| \leq \|u_i - \widetilde{u}_i\| \|u_j\| + \|\widetilde{u}_i\|\|u_j - \widetilde{u}_j\|,$$
and therefore
$$\sum_{i, j \notin A_0} t_i t_j |E_{ij}| \leq \left ( \sum_{i \notin A_0} t_i \|u_i - \widetilde{u}_i\| \right ) \left ( \sum_{j \notin A_0} t_j \|u_j\|\right ) + \left ( \sum_{j \notin A_0} t_j \|u_j - \widetilde{u}_j\| \right )\left ( \sum_{i \notin A_0} t_i \|\widetilde{u}_i\|\right ).$$
By Cauchy--Schwarz and \eqref{conv2} we have $\sum_{i} t_i \|u_i - \widetilde{u}_i\| \to 0$. Similarly, $\sum_{j} t_j \|u_j\| \leq \sqrt{n}$ by the fact that $\sum_{j} \|u_j\|^2=n$. Moreover, we also have $\sum_{i} \|\widetilde{u}_i\|^2 = \sum_{i} nt_i^2 \leq n$, so the same estimate applies to the second term. Thus, we obtained 
\begin{equation}
\label{eqij}
\sum_{i, j \notin A_0} t_i t_j |E_{ij}| \to 0
\end{equation}
and by \eqref{eqdeltas} we also have
$$ \sum_{i,j \notin A_0} t_i t_j |\langle \widetilde{u}_i, \widetilde{u}_j \rangle| |\Delta_{ij}|^2 \to 0.$$
From this the convergence to $0$ of the right-hand side of \eqref{eqbj} follows immediately, since $|\langle \widetilde{u}_i, \widetilde{u}_j \rangle| = n t_it_j | \langle w_{s(i)}, w_{s(j)}\rangle| \geq n \varphi t_it_j$ for all $i, j \notin A_0$. This completes the proof.

\end{proof}

The following lemma may be viewed as a stability version of the equality case in the triangle inequality for scalars.

\begin{lem}
\label{stabtriangle}
Let $y_1, \ldots, y_N \in \mathbb{K}$ and $z_1, \ldots, z_N \in \mathbb{K}$ be such that $|y_i| \leq 1$ for any $i=1, \ldots, N$. Define
$$M = |z_1| + \ldots + |z_N|.$$
If 
$$y_1z_1 + \ldots + y_Nz_N$$
is a non-negative real number greater than or equal to $ M - \varepsilon$ for some $\varepsilon>0$, then 
$$\sum_{i=1}^{N} |z_i| \left |\sgn(z_i) - y_i \right | \leq \sqrt{2 \varepsilon M}.$$
\end{lem}

\begin{proof} From the Cauchy--Schwarz inequality and the assumption it follows that
\begin{align*}\left ( \sum_{i=1}^{N} |z_i| \left |\sgn(z_i) - y_i \right |\right )^2 &\leq \left ( \sum_{i=1}^{N} |z_i| \right ) \left ( \sum_{i=1}^{N} |z_i| \left |\sgn(z_i) - y_i \right |^2 \right ) \\
&=M \left ( \sum_{i=1}^{N} |z_i| \left |\sgn(z_i) - y_i \right |^2 \right ).
\end{align*}
Therefore, it is enough to establish that 
$$\sum_{i=1}^{N} |z_i| \left |\sgn(z_i) - y_i \right |^2 \leq 2 \varepsilon.$$
For every $i=1,\ldots,N$, let $w_i=y_i\sgn(z_i)^{-1}$. Then $ y_i z_i = |z_i|\,w_i$ and hence
$$\left|\sum_{i=1}^{N} y_i z_i\right| = \sum_{i=1}^{N} y_i z_i\  = \sum_{i=1}^{N} |z_i|\,w_i = \left|\sum_{i=1}^{N} |z_i|\,w_i\right|. $$ Moreover $ |w_i|=|y_i|\le 1 $ for every $i=1,\ldots,N$ and
$$\re\left ( \sum_{i=1}^{N}|z_i|\,w_i \right ) = \sum_{i=1}^{N}|z_i|\,w_i \ge M-\varepsilon = \sum_{i=1}^{N}  |z_i| - \varepsilon. $$ 
It follows that 
\begin{equation} 
\label{stabtriangleeq1} 
\sum_{i=1}^{N} |z_i| \bigl( 1-\re(w_i) \bigr) \le \varepsilon. 
\end{equation} 
For every $i$ we have $|w_i|\leq 1$, and therefore
$$ |1-w_i|^2 = 1-2\re(w_i)+|w_i|^2 \leq 2\bigl(1-\re(w_i)\bigr).$$ 
Multiplying by $|z_i|$ and summing over $i$, we obtain
$$\sum_{i=1}^{N} |z_i|\,|1-w_i|^2 \le  2 \sum_{i=1}^{N} |z_i| \bigl( 1-\re(w_i) \bigr).$$
Combining this estimate with \eqref{stabtriangleeq1}, we get 
$$ \sum_{i=1}^{N} |z_i|\,|1-w_i|^2 \le 2\varepsilon. $$ 
Since $ w_i=y_i\sgn(z_i)^{-1}$, we have 
$$ |1-w_i| = \left| 1-\,y_i\sgn(z_i)^{-1} \right| = \left | \sgn(z_i) - y_i \right |. $$ 
Hence
$$ \sum_{i=1}^{N} |z_i| \left | \sgn(z_i) - y_i \right |^2 \le 2\varepsilon. $$ 
This completes the proof. 

\end{proof}

The last auxiliary result needed before the proof of Theorem~\ref{twrchar} in the general case is very simple, but it captures the key asymptotic estimate.

\begin{lem}
\label{lemasymptotics}
Let $N \geq 1$ be an integer and let $t_1, \ldots, t_N$ be non-negative reals with $\sum_{j=1}^{N} t_j^2 = 1$. For $j, k \in \{1, \ldots, N\}$ let $\Delta_{j, k} \in [0, 2]$. Then
$$\sum_{j=1}^{N} t_j^2 \left ( \sum_{k=1}^{N} t^2_k \Delta_{j, k}   \right )^2 \leq 2\sum_{j, k=1}^{N} t^2_jt^2_k \Delta_{j, k}.$$ 
\end{lem}
\begin{proof}
Since for all $j, k$ we have $\Delta_{j, k} \leq 2$ it follows that
$$\sum_{k=1}^{N} t^2_k \Delta_{j, k} \leq 2 \sum_{k=1}^{N}t_k^2 = 2.$$
Therefore
$$\sum_{j=1}^{N} t_j^2 \left ( \sum_{k=1}^{N} t^2_k \Delta_{j, k} \right )^2 \leq 2\sum_{j=1}^{N} t_j^2  \left ( \sum_{k=1}^{N} t^2_k \Delta_{j, k} \right ) = 2 \sum_{j, k=1}^{N} t^2_jt^2_k \Delta_{j, k}.$$
\end{proof}

\begin{proof}[Proof of Theorem \ref{twrchar} in the general case.] We first prove that the stated inclusions imply that $X$ has maximal projection constant, which turns out to be a straightforward consequence of the polyhedral case. Indeed, since the absolute projection constant is invariant under isometry, without loss of generality we may assume that
$$\absconv \{w_1, \ldots, w_d\} \subseteq B_{X^*} \subseteq  \frac{n}{d\delta_{\mathbb{K}}(n)} Z(w_1, \ldots, w_d).$$
Thus, the norm of $X$ can be expressed as
$$\|x\|_{X} =  \sup \left \{|\langle x, w_1 \rangle |, \ldots, |\langle x, w_d \rangle |,  |\langle x, v_1 \rangle|, |\langle x, v_2 \rangle|, \ldots \right \},$$
for some sequence $(v_i)_{i=1}^{\infty} \subseteq \frac{n}{d\delta_{\mathbb{K}}(n)} Z(w_1, \ldots, w_d)$ dense in the unit sphere of $X^{*}$. For each integer $m \geq 1$ let us consider an $n$-dimensional normed space $X_m$ with the norm given as
$$\|x\|_{X_m} =  \max \left \{|\langle x, w_1 \rangle |, \ldots, |\langle x, w_d \rangle |,  |\langle x, v_1 \rangle|, \ldots, |\langle x, v_m \rangle| \right \}.$$
Then $X_m$ is a polyhedral space, which satisfies the given inclusions, and hence, the equality $\lambda(X_m) = \delta_{\mathbb{K}}(n)$ holds by the previously established polyhedral case. Clearly, we have $X_m \to X$ in the sense of the Banach--Mazur distance. Since the absolute projection constant is continuous with respect to $\dbm$, it follows that $\lambda(X)=\delta_{\mathbb{K}}(n)$.

The main difficulty in the general case is proving the necessity of the inclusions. Let $X$ be a general $n$-dimensional normed space over $\mathbb{K}$ satisfying $\lambda(X) = \delta_{\mathbb{K}}(n)$ and let us pick a sequence of $n$-dimensional polyhedral spaces $X_m \subseteq \ell_{\infty}^{N_m}$ (where $N_m \geq n$ is an integer) with $\dbm(X, X_m) \to 1$. Since the absolute projection constant is continuous with respect to $\dbm$, we know that 
$$\delta_m = \lambda(X_m) \to \lambda(X) = \delta_{\mathbb{K}}(n)$$
as $m \to \infty$. We can further assume that $B_{X^*_m} \subseteq B_{X^*} \subseteq (1+\kappa_m)B_{X^*_m}$ for some sequence $\kappa_m \to 0$. We shall prove that there exist a maximal ETF $w_1, \ldots, w_d$ in $\mathbb{K}^n$ associated with a suitable subsequence of the approximating spaces $X_m$ (not relabeled), a corresponding sequence of linear mappings $T_m:\mathbb{K}^n \to \mathbb{K}^n$ and a sequence of non-negative numbers $\tau_m \to 0$ such that the following approximate inclusions hold
\begin{equation}
\label{approxincl}
(1-\tau_m)\absconv \{w_1, \ldots, w_d\} \subseteq T_m(B_{X_m^*}) \subseteq (1+\tau_m) \frac{n}{d\delta_{\mathbb{K}}(n)} Z(w_1, \ldots, w_d).
\end{equation}
In particular, since for sufficiently large $m$ we will then have
$$\frac{1}{2}\absconv \{w_1, \ldots, w_d\} \subseteq T_m(B_{X^*}) \subseteq 2\frac{n}{d\delta_{\mathbb{K}}(n)} Z(w_1, \ldots, w_d),$$
so this will show that the sequence of operators $(T_m)$ satisfies $\sup_{m} \|T_m\| \leq C$ and $\sup_{m} \|T^{-1}_m\| \leq C$. This follows from the fact that there exist $r, R>0$ independent of a particular maximal ETF such that
$$r\mathbb{B}^n \subseteq \absconv \{w_1, \ldots, w_d\} \quad \text{ and } \quad \frac{n}{d\delta_{\mathbb{K}}(n)} Z(w_1, \ldots, w_d) \subseteq R\mathbb{B}^n,$$
where $\mathbb{B}^n$ denotes the Euclidean ball in $\mathbb{K}^n$. Indeed, for the lower inclusion, we may take $r=\frac{1}{\sqrt{n}}$ since the tight frame property of maximal ETF yields
$$\|x\|^2 = \frac{n}{d} \sum_{i=1}^{d} |\langle x, w_i \rangle |^2  \leq n \max_{1 \leq i \leq d}\{ |\langle x, w_i \rangle|^2\}.$$
For the outer inclusion we can take $R=\frac{n}{\delta_{\mathbb{K}}(n)}$, as a simple estimate by the triangle inequality gives us
$$\|a_1w_1 + \ldots + a_dw_d\| \leq |a_1| + \ldots + |a_d| \leq d.$$
for all scalars $a_i$ with $|a_i| \in [0, 1]$. 
Thus, the inclusions \eqref{approxincl} imply the precompactness of the family of operators $(T_m)$ and the uniform upper bound on the norm of inverse operators. Therefore, since the space of linear operators on $\mathbb{K}^{n}$ is finite-dimensional, there exists a subsequence of $(T_m)$ converging to an invertible linear operator $T: \mathbb{K}^n \to \mathbb{K}^n$. Since by \eqref{approxincl} we also have
$$(1-\tau_m)\absconv \{w_1, \ldots, w_d\} \subseteq T_m(B_{X^*}) \subseteq (1+\tau_m)(1+\kappa_m) \frac{n}{d\delta_{\mathbb{K}}(n)} Z(w_1, \ldots, w_d).$$
it will then follow that
$$\absconv \{w_1, \ldots, w_d\} \subseteq T(B_{X^*}) \subseteq \frac{n}{d\delta_{\mathbb{K}}(n)} Z(w_1, \ldots, w_d).$$
Hence, our main goal is to establish the approximate inclusions \eqref{approxincl} hold for a suitable subsequence.

From now on, we regard $X_m$ as a subspace of $\ell_{\infty}^{N_m}$ via the standard isometric embedding. Let us note that the first part of the reasoning in the polyhedral case that concludes with equation \eqref{norma} can be applied for any polyhedral space $X$ satisfying $\lambda(X) > \delta_{\mathbb{K}}(n-1)$. Let $E_m: \ell_{\infty}^{N_m} \to \ell_{\infty}^{N_m}$ be a linear operator like in Lemma \ref{lemcm}, that is $E_m(X_m) \subseteq X_m$, $\sum_{i=1}^{N_m} \|E_m(e_i)\|_{\infty} = 1$ and $\tr(E_m|_{X_m})=\delta_m$. For $1 \leq i \leq N_m$ we denote $t^{(m)}_i = \sqrt{\|E_m(e_i)\|_{\infty}}$ so that
$$\sum_{i=1}^{N_m} (t_i^{(m)})^2 = 1$$ 
and
$$|E_m(e_i)_j| \leq \|E_m(e_i)\|_{\infty} = (t^{(m)}_i)^2 \quad \text { for all } i, j \in \{1, \ldots, N_m\}.$$
Let $S_m \subseteq \{1, 2, \ldots, N_m\}$ be the set of all indices $1 \leq i \leq N_m$ such that $E_m(e_i) \neq 0$ (or equivalently $t^{(m)}_i \neq 0$). Since $\delta_m \to \delta_{\mathbb{K}}(n)$ for sufficiently large $m$ we have $\delta_m> \delta_{\mathbb{K}}(n-1)$ and the previous reasoning for the non-degeneracy of the defined inner product still applies. Hence, similarly as before, we obtain:
\begin{itemize}
\item The formula $\langle x, y \rangle_{E_m} = \sum_{i=1}^{N_m} (t_i^{(m)})^2 x_i \overline{y_i}$ for $x, y \in X_m$ defines a non-degenerate inner product on the space $X_m$ with an orthonormal basis $x_1^{(m)}, \ldots, x_n^{(m)} \in X_m$.
\item We can represent and estimate $\delta_m$ as 
$$\delta_m = \sum_{i=1}^{n} \left \langle E_m(x^{(m)}_i), x^{(m)}_i \right \rangle_{E_m} \leq \sum_{k, j=1}^{N_m} t^{(m)}_kt^{(m)}_j | \langle u_k^{(m)}, u_j^{(m)} \rangle |,$$
where vectors $u^{(m)}_1, \ldots, u^{(m)}_{N_m} \in \mathbb{K}^n$ are defined by $u^{(m)}_{ij}=t^{(m)}_i \overline{x^{(m)}_{ji}}$ for $1 \leq i \leq N_m$ and $1 \leq j \leq n$.
\item For every vector $\sum_{i=1}^{n} c_ix^{(m)}_i \in X_m$ we have 
$$\left \| \sum_{i=1}^{n} c_ix^{(m)}_i \right \|_{\infty} = \max_{1 \leq j \leq N_m} |\langle c, f^{(m)}_j \rangle|,$$
where $c=(c_1, \ldots, c_n) \in \mathbb{K}^n$ and the functionals $f^{(m)}_j \in \mathbb{K}^n$ are given as $f^{(m)}_j = \left ( \overline{x^{(m)}_{1j}}, \ldots, \overline{x^{(m)}_{nj}} \right )$ for $j=1, \ldots, N_m$.
\end{itemize}
We note that by the last formula, the basis $x_1^{(m)}, \ldots, x_n^{(m)}$ induces the linear isomorphism $T_m$ required in \eqref{approxincl}, as $T_m({B_{X^*_m}}) = \absconv\{f_1^{(m)},\ldots,f_{N_m}^{(m)}\}$.

After passing to an appropriate subsequence if necessary, we may assume that all conclusions of Theorem~\ref{twrboundstab} hold. Therefore, we are given the partition  $\{1, \ldots, N_m\} = A^{(m)}_0 \sqcup A^{(m)}_1 \sqcup \ldots \sqcup A^{(m)}_d$ of all the indices. Since $t^{(m)}_i=0$ implies $u^{(m)}_i=0$ by definition, it follows that $\{1, \ldots, N_m\} \setminus S_m \subseteq A^{(m)}_0$, since the exceptional set $A^{(m)}_0$ contains all zero vectors $u^{(m)}_i$. 

We start by establishing the approximate lower inclusion. Let us fix $s \in \{1, \ldots, d\}$. By condition $(3)$ of Theorem~\ref{twrboundstab} we have 
\[ \sum_{k\in A_s^{(m)}} (t_k^{(m)})^2 \to \frac1d. \]
Moreover, conditions $(2)$ and $(4)$ yield together
$$ \sum_{k\in A_s^{(m)}} (t_k^{(m)})^2 \left ( (d_k^{(m)})^2 + (v_k^{(m)})^2 \right ) \to 0,$$
where 
$$d^{(m)}_k = \inf_{|\alpha|=1} \left\| \frac{u_k^{(m)}}{\|u_k^{(m)}\|} -\alpha w_s \right \| \quad \text{ and } \quad v_k^{(m)}= \left | \frac{\|u^{(m)}_k\|}{t^{(m)}_k} - \sqrt{n} \right |.$$ 
Let $j=j_s^{(m)} \in A_s^{(m)}$ be an index for which
$$ (d_j^{(m)})^2 + (v_j^{(m)})^2 = \min_{k \in A_s^{(m)}} \left ( (d_k^{(m)})^2 + (v_k^{(m)})^2 \right ).$$
Then
\begin{align*}
 \sum_{k\in A_s^{(m)}} (t_k^{(m)})^2 \left ( (d_k^{(m)})^2 + (v_k^{(m)})^2 \right ) &\geq  \sum_{k\in A_s^{(m)}} (t_k^{(m)})^2 \left ( (d_j^{(m)})^2 + (v_j^{(m)})^2 \right ) \\
&\geq \frac{1}{2d}\left ( (d_j^{(m)})^2 + (v_j^{(m)})^2 \right ),
\end{align*}
since $\sum_{k\in A_s^{(m)}} (t_k^{(m)})^2 \geq \frac{1}{2d}$ for sufficiently large $m$. This shows that, for the chosen index $d_{j^{(m)}_s}^{(m)} \to 0$ and $v_{j^{(m)}_s}^{(m)} \to 0$. Moreover, passing to a subsequence if necessary, we can further assume that the sequence of unimodular scalars realizing the infimum in $d_{j^{(m)}_s}^{(m)}$ converges to a unimodular scalar $\alpha_s$. Since 
\[\frac{u_{j_s^{(m)}}^{(m)}} {\left \| u_{j_s^{(m)}}^{(m)} \right \|} \to \alpha_sw_s \quad \text{ and } \quad  \frac{\left \|u_{j_s^{(m)}}^{(m)} \right \|}{t_{j_s^{(m)}}^{(m)}} \to \sqrt n, \]
we get 
\[  \frac{u_{j_s^{(m)}}^{(m)}} {t_{j_s^{(m)}}^{(m)}} \to \sqrt n\,\alpha_s w_s.\]
However, as the functional $f^{(m)}_k$ is given as
\[ f_{k}^{(m)} = \frac{u_k^{(m)}}{t_k^{(m)}}, \] 
for $k \in A_s^{(m)}$, we obtain that for every $c=(c_1,\dots,c_n)\in\mathbb K^n$ 
\[\left | \left \langle c, f_{j_s^{(m)}}^{(m)} \right \rangle \right | = \left | \left\langle c, \frac{u_{j_s^{(m)}}^{(m)}} {t_{j_s^{(m)}}^{(m)}} \right\rangle \right | \to  \sqrt n\,|\langle c,w_s\rangle|. \]
Thus, up to the factor $\sqrt n$, for every $s=1,\ldots,d$ each of the functionals $x\mapsto |\langle x,w_s\rangle|$ appears in the limit (after passing to further subsequences for each $s$, if necessary). Hence, since the absolutely convex hull depends continuously on finitely many vertices, there exists a sequence of non-negative numbers $\tau_m \to 0$ such that $(1-\tau_m)\absconv \{w_1, w_2, \ldots, w_d\} \subseteq T_m(B_{X_m^*})$ for a suitable linear transformation $T_m$. This proves the approximate lower inclusion.

It remains to prove the approximate outer zonotope inclusion. The key step is to establish the following convergence:
\begin{equation}
\label{approxid}
\|E_m(x^{(m)}_i) - Cx^{(m)}_i\|_{E_m} \to 0
\end{equation}
for any $i=1, \ldots, n$. To this end, for $j, k \in S_m$ we define 
$$y_{jk} = \frac{E_m(e_k)_j}{(t^{(m)}_k)^2} \qquad \text{ and } \qquad z_{jk} = t^{(m)}_jt^{(m)}_k \langle u^{(m)}_j, u^{(m)}_k \rangle.$$
Then $|y_{jk}| \leq 1$ for all $j, k \in S_m$ and
\begin{align*}
\delta_m &= \sum_{i=1}^{n} \langle E_m(x^{(m)}_i), x^{(m)}_i \rangle_{E_m} = \sum_{j, k=1}^{N_m} \sum_{i=1}^{n} (t^{(m)}_j)^2E_m(e_k)_j x^{(m)}_{ik} \overline{x^{(m)}_{ij}} \\ 
&= \sum_{j, k \in S_m} \sum_{i=1}^{n} (t^{(m)}_j)^2E_m(e_k)_j x^{(m)}_{ik} \overline{x^{(m)}_{ij}}=\sum_{j, k \in S_m}  \frac{\langle u^{(m)}_j, u^{(m)}_k \rangle}{t^{(m)}_jt^{(m)}_k}(t^{(m)}_j)^2E_m(e_k)_j \\
&= \sum_{j, k \in S_m} \frac{E_m(e_k)_j}{(t^{(m)}_k)^2} t^{(m)}_jt^{(m)}_k \langle u^{(m)}_j, u^{(m)}_k \rangle = \sum_{j, k \in S_m} y_{jk} z_{jk}.
\end{align*}
Since $\sum_{j, k \in S_m} |z_{jk}| \leq \delta_{\mathbb{K}}(n)$, it follows that 
$$\sum_{j, k \in S_m} y_{jk} z_{jk}  \geq \sum_{j, k \in S_m} |z_{jk}| - \varepsilon_m,$$
where $\varepsilon_m = \delta_{\mathbb{K}}(n) - \delta_m$. Hence, Lemma \ref{stabtriangle} yields 
$$\sum_{j, k \in S_m}  t^{(m)}_jt^{(m)}_k |\langle u^{(m)}_j, u^{(m)}_k \rangle| \Delta^{(m)}_{j, k} \leq \sqrt{2 \varepsilon_m \sum_{j, k \in S_m} t^{(m)}_jt^{(m)}_k |\langle u^{(m)}_j, u^{(m)}_k \rangle| } \leq \sqrt{2\delta_{\mathbb{K}}(n) \varepsilon_m},$$
where 
$$\Delta^{(m)}_{j, k} = \left |\frac{E_m(e_k)_j}{(t^{(m)}_k)^2} - \sgn \langle u^{(m)}_j, u^{(m)}_k \rangle \right |,$$
for $j, k \in S_m$. Clearly, we have $\Delta^{(m)}_{j, k} \leq 2$ for all $j, k$ .

For $s, r \in \{1, \ldots, d\}$ and $j \in A_s^{(m)}$, $k \in A_r^{(m)}$ let us now write
$$\langle u_j^{(m)},u_k^{(m)}\rangle = n\,t_j^{(m)}t_k^{(m)} \alpha_j^{(m)} \overline{\alpha_k^{(m)}} \langle w_s,w_r\rangle + E_{jk}^{(m)}.$$
By estimate \eqref{eqij} from Lemma \ref{lemtiuistab} we then have
$$\sum_{j, k \notin A^{(m)}_0} t^{(m)}_jt^{(m)}_k |E_{jk}^{(m)}| \to 0.$$
Therefore, since $ |\langle w_r,w_s\rangle| \geq \varphi_{\mathbb{K}}(n) >0$, we now have
$$\sum_{j, k \notin A^{(m)}_0} t^{(m)}_jt^{(m)}_k |\langle u^{(m)}_j, u^{(m)}_k \rangle| \Delta^{(m)}_{j, k}$$
$$\geq n \varphi_{\mathbb{K}}(n) \sum_{j, k \notin A^{(m)}_0} (t_j^{(m)})^2(t_k^{(m)})^2 \Delta^{(m)}_{j, k} - 2\sum_{j, k \notin A^{(m)}_0} t^{(m)}_jt^{(m)}_k|E_{jk}^{(m)}|,$$
which, together with the previous convergences, implies
$$\sum_{j, k \notin A_0^{(m)}} (t_j^{(m)})^2(t_k^{(m)})^2 \Delta^{(m)}_{j, k} \to 0.$$
Since $\Delta^{(m)}_{j, k} \leq 2$ for all $j, k \in S_m$ and $\sum_{j \in A_0^{(m)}} (t_j^{(m)})^2 \to 0$, the same convergence remains true when the indices from $S_m \cap A_0^{(m)}$ are included. Hence we conclude
\begin{equation}
\label{eqdeltaconv}
\sum_{j, k \in S_m} (t_j^{(m)})^2(t_k^{(m)})^2 \Delta^{(m)}_{j, k} \to 0.
\end{equation}

Let us now fix $i \in \{1, \ldots, n\}$ and calculate
\begin{align*}
\|E_m(x_i^{(m)})-Cx_i^{(m)}\|^2_{E_m} &=  \sum_{j=1}^{N_m} (t_j^{(m)})^2 \left| E_m(x_i^{(m)})_j - C\,x_{ij}^{(m)} \right|^2 \\
 &=\sum_{j \in S_m} (t_j^{(m)})^2 \left| E_m(x_i^{(m)})_j - C\,x_{ij}^{(m)} \right|^2.
\end{align*}
 For a given $j \in S_m$ we have
\begin{align*}
E_m(x_i^{(m)})_j - C\,x_{ij}^{(m)} &= \sum_{k=1}^{N_m}E_m(e_k)_jx_{ik}^{(m)} - Cx_{ij}^{(m)}=\sum_{k \in S_m}E_m(e_k)_jx_{ik}^{(m)} - Cx_{ij}^{(m)} \\
 &=\sum_{k \in S_m}\frac{E_m(e_k)_j}{t^{(m)}_k} \overline{u_{ki}^{(m)}} - \frac{C}{t^{(m)}_j}{\overline{u_{ji}^{(m)}}}
\end{align*}
 and therefore
 $$\|E_m(x_i^{(m)})-Cx_i^{(m)}\|^2_{E_m} = \sum_{j \in S_m} \left | \sum_{k \in S_m} t^{(m)}_j\frac{E_m(e_k)_j}{t^{(m)}_k} \overline{u_{ki}^{(m)}} - C{\overline{u_{ji}^{(m)}}} \right |^2 $$
 $$\leq 2\sum_{j \in S_m} \left |t^{(m)}_j\sum_{k \in S_m} \frac{E_m(e_k)_j}{t^{(m)}_k} \overline{u_{ki}^{(m)}} -  t_j^{(m)} \sum_{k \in S_m} t_k^{(m)} \overline{ \sgn \langle u_k^{(m)}, u_j^{(m)}\rangle u_{ki}^{(m)}}  \right |^2$$
 $$+2\sum_{j \in S_m}\left |C{\overline{u_{ji}^{(m)}}}  -  t_j^{(m)} \sum_{k \in S_m} t_k^{(m)}\overline{\sgn \langle u_k^{(m)}, u_j^{(m)}\rangle u_{ki}^{(m)}}  \right |^2.$$
The second term is obviously bounded above by
$$ \sum_{j \in S_m} \left \|C{u_{j}^{(m)}}  -  t_j^{(m)} \sum_{k \in S_m} t_k^{(m)} \sgn \langle u_k^{(m)}, u_j^{(m)}\rangle u_{k}^{(m)} \right \|^2,$$
which tends to $0$ by Lemma \ref{lemtiuistab}. To establish \eqref{approxid} it is therefore enough to show the convergence of the first term. We can estimate it as
$$\sum_{j \in S_m} (t^{(m)}_j)^2 \left | \sum_{k \in S_m} \left ( \frac{E_m(e_k)_j}{t^{(m)}_k} - t_k^{(m)} \sgn \langle u_k^{(m)}, u_j^{(m)} \right )u^{(m)}_{ki} \right |^2$$
$$\leq \sum_{j \in S_m} (t^{(m)}_j)^2 \left | \sum_{k \in S_m} t^{(m)}_k \Delta_{j, k}|u^{(m)}_{ki}| \right |^2 \leq \sum_{j \in S_m} (t^{(m)}_j)^2 \left | \sum_{k \in S_m} t^{(m)}_k \Delta^{(m)}_{j, k}\|u_k^{(m)}\| \right |^2.$$
For $k \in S_m$ let us write $\|u^{(m)}_k\|=\sqrt{n}t^{(m)}_k + r^{(m)}_k$. Then
$$\left | \sum_{k \in S_m} t^{(m)}_k \Delta^{(m)}_{j, k}\|u^{(m)}_k\| \right |^2 \leq 2n \left | \sum_{k \in S_m} (t^{(m)}_k)^2 \Delta^{(m)}_{j, k} \right |^2 + 2 \left | \sum_{k \in S_m} t^{(m)}_k \Delta^{(m)}_{j, k}r^{(m)}_k \right |^2.$$
For the first term, by \eqref{eqdeltaconv} and Lemma \ref{lemasymptotics} we have
$$\sum_{j \in S_m} (t^{(m)}_j)^2 \left | \sum_{k \in S_m} (t^{(m)}_k)^2 \Delta^{(m)}_{j, k} \right |^2 \to 0. $$
For the second term, it is enough to use the estimate $\Delta^{(m)}_{j, k} \leq 2$, the Cauchy--Schwarz inequality and condition $(4)$ of Theorem \ref{twrboundstab}, as
$$\sum_{j \in S_m} (t^{(m)}_j)^2 \left | \sum_{k \in S_m} t^{(m)}_k \Delta^{(m)}_{j, k}r^{(m)}_k \right |^2 \leq 4\sum_{j \in S_m} (t^{(m)}_j)^2 \left | \sum_{k \in S_m} t^{(m)}_k |r^{(m)}_k| \right |^2$$
$$\leq 4 \left ( \sum_{j \in S_m} (t^{(m)}_j)^2 \right ) \left ( \sum_{k \in S_m} (t_k^{(m)})^2 \right ) \left ( \sum_{k \in S_m} |r^{(m)}_k|^2 \right ) \leq 4 \sum_{k \in S_m} |r^{(m)}_k|^2 \to 0.$$
This concludes the proof of the convergence \eqref{approxid}.

Now we shall deduce the approximate outer zonotope inclusion from the convergence \eqref{approxid}. Since both inclusions are rescaled by $\sqrt n$, it is enough to establish that for every $j \in \{1, \ldots, N_m\}$ we have
$$f^{(m)}_j = \left ( \overline{x^{(m)}_{1j}}, \ldots, \overline{x^{(m)}_{nj}} \right ) \in \sqrt{n}(1+\tau_m) \frac{n}{d\delta_{\mathbb{K}}(n)} Z(w_1, \ldots, w_d).$$
where $\tau_m \to 0$. For a fixed $j \in \{1, \ldots, N_m\}$ we can write
$$(E_m(x_1^{(m)})_{j},\ldots,E_m(x_n^{(m)})_{j}) = \sum_{k \in S_m} E_m(e_k)_{j} (x_{1k}^{(m)},\dots,x_{nk}^{(m)}).$$
By Theorem \eqref{twrboundstab}, for every $s \in \{1, \ldots, d\}$ and $k \in A_s^{(m)}$ we have 
\[ (x_{1k}^{(m)},\dots,x_{nk}^{(m)}) = \frac{\overline{u_k^{(m)}}}{t_k^{(m)}} = \sqrt n\, \overline{\alpha_k^{(m)} w_s} + r_k^{(m)}, \] 
where $|\alpha_k^{(m)}|=1$ and 
\[ \sum_{k\in A_s^{(m)}} (t_k^{(m)})^2 \|r_k^{(m)}\|^2 \to 0. \] 
Therefore 
\[ (E_m(x_1^{(m)})_{j},\dots,E_m(x_n^{(m)})_{j}) = \sqrt n \sum_{s=1}^{d} b_s^{(m)} \overline{w_s} + \sum_{s=0}^{d} v^{(m)}_s, \] where 
\[ b_s^{(m)} = \sum_{k\in A_s^{(m)}} \overline{\alpha_k^{(m)}} E_m(e_k)_{j} \quad \text{ and } \quad v^{(m)}_s = \sum_{k\in A_s^{(m)}} E_m(e_k)_{j} r_k^{(m)}\]
for $s \in \{1, \ldots, d\}$ and we additionally define $v_0^{(m)}$ as
$$v_0^{(m)} = \sum_{k \in A_0^{(m)} \cap S_m} E_m(e_k)_{j} ( x_{1k}^{(m)}, \dots, x_{nk}^{(m)} ).$$
Using $|E_m(e_k)_{j}| \le (t_k^{(m)})^2$, we get 
$|b_s^{(m)}| \le \sum_{k\in A_s^{(m)}} (t_k^{(m)})^2$
and by the triangle inequality and the Cauchy--Schwarz inequality
$$\|v^{(m)}_s\| \leq \sum_{k\in A_s^{(m)}} (t_k^{(m)})^2 \|r_k^{(m)}\| \leq \left ( \sum_{k\in A_s^{(m)}} (t_k^{(m)})^2 \right )^{\frac{1}{2}} \left ( \sum_{k\in A_s^{(m)}} (t_k^{(m)})^2 \|r_k^{(m)}\|^2 \right )^{\frac{1}{2}} \to 0.$$
Moreover,
$$\|v_0^{(m)}\| = \left \| \sum_{k \in A_0^{(m)} \cap S_m} \frac{E_m(e_k)_{j}}{t_k^{(m)}}\overline{u_k^{(m)}} \right \| \leq \sum_{k \in A_0^{(m)} \cap S_m} t_k^{(m)} \|u_k^{(m)}\| \to 0,$$
by the Cauchy--Schwarz inequality and condition $(1)$ of Theorem \ref{twrboundstab}. Therefore, there exists a sequence $\gamma_m \to 0$ such that 
$$\sum_{s=0}^{d} v^{(m)}_s \in \gamma_m \frac{\sqrt{n}}{d}  Z(\overline{w}_1, \ldots, \overline{w}_d).$$ 
Furthermore, by condition $(3)$ of Theorem~\ref{twrboundstab} we can estimate $|b_s^{(m)}| \leq \frac{1}{d} + \omega_m$ for some sequence $\omega_m \to 0$. Therefore, for any $j \in \{1, \ldots, N_m\}$ we have
$$ \left( E_m(x_1^{(m)})_{j},\dots,E_m(x_n^{(m)})_{j} \right ) \in \frac{\sqrt{n}}{d}(1+\omega_m+\gamma_m) Z(\overline{w}_1, \ldots, \overline{w}_d),$$
where $\omega_m$ and $\gamma_m$ do not depend on the index $j$. Now, let us write $E_m(x_i^{(m)}) = a^{(m)}_{i1}x_1^{(m)} + \ldots + a^{(m)}_{in}x_n^{(m)}$ and let $A_m=(a^{(m)}_{ij})$ be the corresponding $n \times n$ matrix. Since the vectors $x_1^{(m)}, \ldots, x_n^{(m)}$ form an $E_m$-orthonormal basis, from the established convergence we have
$$\sum_{j=1}^{n} \left | a_{ij}^{(m)} - C \delta_{ij} \right |^2 = \left \| E_m(x_{i}^{(m)}) - Cx_{i}^{(m)}  \right \|^2_{E_m} \to 0,$$
and hence $a^{(m)}_{ij} \to C \delta_{ij}$ (the Kronecker delta) for all $i, j \in \{1, \ldots, n\}$. Thus, the matrix $A_m$ converges entrywise to $C\id_n$ and for sufficiently large $m$ it is invertible with $A_m^{-1} \to C^{-1} \id_n$. Hence, also $(\overline{A_m})^{-1} \to C^{-1} \id_n$. Let us write $(\overline{A_m})^{-1} =C^{-1}\id_n + J_m$, where $\|J_m\| \to 0$. Then
$$(\overline{A_m})^{-1} \left (Z(w_1, \ldots,  w_d)\right ) \subseteq C^{-1}Z(w_1, \ldots w_d) + \rho_m Z(w_1, \ldots w_d)$$
where $\rho_m=\frac{R}{r} \|J_m\|$ and $0<r<R$ are such that 
$$r\mathbb{B}^n \subseteq Z(w_1, \ldots,  w_d) \subseteq R \mathbb{B}^n.$$
Therefore, since
$$\overline{A_m}\left ( \overline{x^{(m)}_{1j}}, \ldots, \overline{x^{(m)}_{nj}} \right ) = \left (\overline{E_m(x_1^{(m)})_{j}}, \ldots, \overline{E_m(x_n^{(m)})_{j}} \right )\in  \frac{\sqrt{n}}{d}(1+\omega_m+\gamma_m) Z(w_1,\dots,w_d)$$
for every $j=1, \ldots, N_m$ and sufficiently large $m$ for the functional $f^{(m)}_j = \left ( \overline{x^{(m)}_{1j}}, \ldots, \overline{x^{(m)}_{nj}} \right )$ we have
\begin{align*}
f^{(m)}_j &\in (\overline{A_m})^{-1} \left ( \frac{\sqrt{n}}{d}(1+\omega_m+\gamma_m) Z(w_1,\dots,w_d) \right ) \\
&\subseteq C^{-1}\frac{\sqrt{n}}{d}(1+C\rho_m)(1+\omega_m+\gamma_m)Z(w_1,\dots,w_d) \\
&=\sqrt{n}\frac{n}{d\delta_{\mathbb{K}}(n)}(1+\tau_m) Z(w_1, \ldots, w_d).
\end{align*}
This concludes the proof of the characterization in the general case.

\end{proof}

\section{Uniqueness of spaces with the maximal projection constant}

Using the characterization given in Theorem~\ref{twrchar} and Lemma~\ref{lemincl}, it is straightforward to deduce that the case $\mathbb{R}^2$ is the only situation considered in this paper in which there is a unique space with the maximal projection constant. K\"onig and Tomczak-Jaegermann attempted to establish an even more general result in Theorem~3 of \cite{konigtomczak2}. Besides the fact that its proof relied on the incorrect Proposition~5, the approach is not particularly effective from our point of view. Rather than using the inclusion~\eqref{inclusion} directly, the authors employed instead an involved probabilistic argument based on estimates of Euclidean norms. Consequently, the result misses many low-dimensional cases, since it requires an embedding into $\ell_\infty^N$ with $8 \leq N \leq e^{\frac{\sqrt{n}}{8e}}$. In particular, for the result to provide any information, one must have $e^{\frac{\sqrt{n}}{8e}} \geq 8$ and none of the currently known real or complex dimensions admitting a maximal ETF satisfy this condition. It also appears that no proof of the uniqueness of the regular hexagon norm in $\mathbb{R}^2$ was provided in that paper. In the following we adopt a more direct approach to the uniqueness question in those dimensions for which a maximal ETF exists.

We first prove the following lemma, which, starting from a pair of non-isometric normed spaces, produces an uncountable family of pairwise non-isometric spaces interpolating between them.

\begin{lem}
\label{lemisometric}
Let $n \geq 2$ and let $X_0 = (\mathbb{K}^n, \| \cdot\|_0)$, $X_1=(\mathbb{K}^n, \| \cdot\|_1)$ be two normed spaces. For $t \in [0, 1]$ we define a normed space $X_t = (\mathbb{K}^n, \| \cdot \|_t)$ with the norm
$$\|x\|_t = (1-t)\|x\|_0 + t\|x\|_1 \qquad \text{ for } x \in \mathbb{K}^n.$$
If the spaces $X_0$ and $X_1$ are non-isometric, then among the spaces $X_t$ for $t\in[0,1]$ there exists a subset of cardinality continuum consisting of pairwise non-isometric spaces.
\end{lem}

\begin{proof}

We define a mapping $f:[0, 1] \to \mathbb{R}$ as
$$f(t) = \dbm(X_t, X_0).$$
Clearly, $f(0)=1$, $f(t) \geq 1$ for every $t \in [0, 1]$ and $f(1)>1$ by the assumption. We claim that $f$ is a continuous function of $t$. If $0 < t < s < 1$, then by the triangle inequality
$$f(s) = \dbm(X_s, X_0) \leq \dbm(X_s, X_t) \cdot \dbm(X_t, X_0) = \dbm(X_s, X_t) \cdot f(t).$$
Similarly $f(t) \leq \dbm(X_s, X_t) \cdot f(s)$. Moreover for any $x \in \mathbb{K}^n$
$$\frac{t}{s} \left( (1-s)\|x\|_0 + s\|x\|_1  \right ) \leq (1-t)\|x\|_0 + t \|x\|_1 \leq \frac{1-t}{1-s} \left ((1-s)\|x\|_0 + s \|x\|_1 \right ).$$
This yields
$$\dbm(X_s, X_t) \leq \frac{s-ts}{t-ts}$$
and consequently
$$1 \leq \max \left \{ \frac{f(t)}{f(s)}, \frac{f(s)}{f(t)} \right \} \leq \frac{s-ts}{t-ts}.$$
Since $\frac{s-ts}{t-ts} \to 1$ for $s \to t$, this proves continuity of $f$ in every point of the open interval $(0, 1)$.

To establish the continuity in the endpoints, let us pick numbers $\beta > \alpha >0$ such that 
$$\alpha \|x\|_1 \leq \|x\|_0 \leq \beta\|x\|_1,$$
for any $x \in \mathbb{K}^n$. Then for $t \in [0, 1]$ and $x \in \mathbb{K}^n$ we have
\begin{align*}
\left ( 1 + \frac{t}{\beta} - t \right ) \|x\|_0 &=  (1-t) \|x\|_0 + \frac{t}{\beta} \|x\|_0 \leq (1-t) \|x\|_0 + t \|x\|_1 = \|x\|_t \\
 &\leq (1-t) \|x\|_0 + \frac{t}{\alpha} \|x\|_0 = \left ( 1 + \frac{t}{\alpha} - t \right ) \|x\|_0.
\end{align*}
Thus $f(t) \leq \frac{ 1 + \frac{t}{\alpha} - t}{1 + \frac{t}{\beta} - t}$. This expression clearly approaches $1$ as $t \to 0^{+}$, which proves the continuity in $f$ also in $0$. For the continuity in $1$, we note that for every $t \in [0, 1]$ and $x \in \mathbb{K}^n$
\begin{align*}
\left ((1-t)\alpha + t \right ) \|x\|_1  &= (1-t)\alpha \|x\|_1 + t \|x\|_1 \leq (1-t)\|x\|_0 + t\|x\|_1 = \|x\|_t \\
&\leq (1-t) \beta \|x\|_1 + t \|x\|_1 = ((1-t)\beta + t) \|x\|_1.
\end{align*}
Hence
$$\dbm(X_t, X_1) \leq \frac{(1-t)\beta + t}{(1-t)\alpha + t}$$
and from the triangle inequality for $\dbm$
$$\frac{\dbm(X_0, X_1)}{\dbm(X_t, X_1)} \leq \dbm(X_t, X_0) = f(t) \leq \dbm(X_0, X_1) \cdot \dbm(X_t, X_1),$$
so that
$$\dbm(X_0, X_1) \cdot \frac{(1-t)\alpha + t}{(1-t)\beta + t} \leq f(t) \leq \dbm(X_0, X_1) \cdot \frac{(1-t)\beta + t}{(1-t)\alpha + t}.$$
Since 
$$\lim_{t \to 1^{-1}} \frac{(1-t)\beta + t}{(1-t)\alpha + t} = 1, $$
this yields the continuity of $f$ also in $1$.

To conclude the proof, observe that $f(0)=1$ and $f(1)>1$ by the assumption. As $f$ is continuous and nonconstant, the image $f([0, 1])$ contains some nondegenerate interval, and therefore, it contains continuum many different values. Since the Banach--Mazur distance to $X_0$ is an isometric invariant, spaces $X_t$ corresponding to different values of $t$ are non-isometric to each other. Hence there are continuum many pairwise non-isometric spaces among the spaces $X_t$. This concludes the proof.

\end{proof}

\begin{cor}
\label{coruniq}
Let $n \geq 2$ be an integer such that there exists a maximal ETF in $\mathbb{K}^n$ (of cardinality $d = d_{\mathbb{K}}(n)$). Then there exists a unique $n$-dimensional normed space $X$ (up to isometry) satisfying $\lambda(X)=\lambda_{\mathbb{K}}(n)$ if and only if $n=2$ and $\mathbb{K}=\mathbb{R}$. In this case, the unit ball of $X$ is an affine image of a regular hexagon. In all other situations, there exist continuum many pairwise non-isometric spaces with the maximal projection constant.
\end{cor}

\begin{proof}
We start with the case $n=2$ and $\mathbb{K}=\mathbb{R}$. It follows immediately from Theorem \ref{twrchar} and Lemma \ref{lemincl} that in this situation, any normed space $X$ with the maximal projection constant has the dual unit ball equal to some $0$-symmetric affine regular hexagon. It is straightforward to check that in this case, the unit ball of $X$ is also a $0$-symmetric affine regular hexagon. Thus, any two-dimensional real spaces with maximal projection constants are isometric.

In all other situations, it follows from Lemma \ref{lemincl} that the inclusion $K \subseteq L$ is strict, where
$$K=\absconv \{w_1, w_2, \ldots, w_d\} \quad \text{ and } \quad L= \frac{n}{d\delta_{\mathbb{K}}(n)} Z(w_1, \ldots, w_d).$$
Let $X_0$ and $X_1$ be normed spaces with the unit balls $K$ and $L$, respectively. We claim that $X_0$ and $X_1$ are non-isometric to each other. Indeed, let us first note that every point $w_j$ for $j=1, \ldots, d$ is an extreme point of $L$. This follows from the fact that if $j \in \{1, \ldots, d\}$ and $x= \frac{n}{d\delta_{\mathbb{K}}(n)} (a_1w_1 + \ldots + a_dw_d) \in L$, where $|a_i| \leq 1$, then
$$|\langle x, w_j \rangle | \leq \frac{n}{d\delta_{\mathbb{K}}(n)} (1 + (d-1)\varphi_{\mathbb{K}}(n)) =  1, $$
with the equality if and only if there exists a unimodular scalar $t$ with $a_j=t$ and $a_i = \sgn \langle w_i, w_j \rangle t$ for all $i \neq j$. However, by Lemma \ref{lemwj} this means that $x=tw_j$. Thus, up to a multiplication by a unimodular scalar, the functional $x \mapsto \langle x, w_j \rangle $ is maximized on $L$ exactly on $w_j$. This means that this functional exposes $w_j$ in $L$, so $w_j$ is indeed an extreme point of $L$.

If $X_1$ were isometric to $X_0$, its unit ball $L$ would be an absolutely convex hull of $d$ points. However, since all the points $w_j$ are extreme, $L$ would have to be equal to exactly $\absconv \{w_1, \ldots, w_d\} = K$. By Lemma \ref{lemincl} we already know that this is impossible, since $K$ is a proper subset of $L$.

For $t \in [0, 1]$ the unit ball of a normed space $X_t$ defined like in Lemma \ref{lemisometric} is between $K$ and $L$, and hence, the dual space $X^*_t$ has the maximal projection constant by Theorem \ref{twrchar}. Therefore, the desired conclusion follows directly from Lemma \ref{lemisometric}, since the dual spaces of non-isometric spaces are non-isometric as well. 
\end{proof}

It should be noted that non-uniqueness of a space the with maximal projection constant may also arise from non-uniqueness of maximal ETFs in $\mathbb K^n$. This occurs, for example, in the case of $\mathbb C^3$. While maximal ETFs are known to be unique up to orthogonal transformations in $\mathbb R^n$ for $n=2,3,7,23$, and unique up to unitary transformations in $\mathbb C^2$, the situation in $\mathbb C^3$ is quite different: there exists an infinite family of pairwise non-equivalent maximal ETFs. For a complete description of all maximal ETFs in $\mathbb C^3$, see Theorem~5.5 of \cite{hughstonsalamon}. The reader is referred to \cite{fickusmixon} for tables summarizing the currently known existence results for equiangular tight frames.

\section*{Acknowledgements} 
The author would like to express sincere gratitude to the anonymous referee for identifying a serious gap in the first version of the manuscript, for providing the example presented in Remark~\ref{remarkcounterexample}, and for numerous constructive suggestions that substantially improved the paper. The author also thanks Microsoft Copilot for assistance during the revision of the manuscript and Florian Grundbacher for help with verifying several arguments in the revised version.

\end{document}